\documentclass[10pt]{article}
\date{}
\usepackage{amssymb}
\usepackage{amsmath}
\usepackage{amsthm}
\usepackage{graphicx}
\usepackage{indentfirst}
\allowdisplaybreaks[4]

\numberwithin{equation}{section}
\begin{document}

\title{Gradient Estimates For The CR Heat Equation On Complete Noncompact Pseudohermitian Manifolds\footnote{
		This work is supported by NSFC Grants No.12171091.}  }
		\author{Yuxin Dong, Yibin Ren and Biqiang Zhao\footnote{
		E-mail addresses: 17110840003@fudan.edu.cn}}
		
\maketitle
\setlength{\parindent}{2em}

\begin{abstract}
In this paper, we derive local and global Li-Yau type gradient estimates for the positive solutions of the CR heat equation on complete noncompact pseudo-Hermitian manifolds. As applications of the gradient estimates, we give a Harnack inequality for the positive solutions of the CR heat equation, and then obtain an upper bound estimate for the corresponding heat kernel.
\end{abstract}


\section{Introduction}
\label{1}
  The method of gradient estimates is an important tool in geometric analysis,
which is originated first in Yau \cite{Ya} and Cheng and Yau \cite{ChY} for investigating harmonic functions and further developed in Li and Yau \cite{LY}
for studying solutions of the heat equation on complete Riemannian
manifolds. In \cite{LY}, Li and Yau established their celebrated parabolic
inequality, which asserts that, if $M$ is an $n$-dimensional complete
Riemannian manifold with Ricci curvature bounded from below by $-K$, for
some constant $K\geq 0$ and $u$ is any positive solution of the heat equation%
\[
(\bigtriangleup -\frac{\partial }{\partial t})u=0 ,
\]%
then%
\[
\frac{|\nabla u|^{2}}{u^{2}}-\alpha \frac{u_{t}}{u}\leq \frac{n\alpha ^{2}}{%
2t}+\frac{n\alpha ^{2}}{2(\alpha -1)}K 
\]%
for all $t>0$, $\alpha >1$. In particular, if $K=0$, then the following more
precise inequality holds%
\[
\frac{|\nabla u|^{2}}{u^{2}}-\frac{u_{t}}{u}\leq \frac{n}{2t}. 
\]%
Since then many improvements or generalizations of Li-Yau's parabolic
inequality have been developed on Riemannian manifolds or more general
metric measure spaces, see e.g. \cite{Da}, \cite{Ha}, \cite{BQ}, \cite{CaN}, \cite{BL}, \cite{SZ}, \cite{LX}, \cite{BBG},
\cite{ZZ}, \cite{CZ}, \cite{ZhZ}, \cite{YZ} and the references therein.
\par
The Li-Yau's inequality has also been generalized for non elliptic operators
that include subelliptic operators on sub-Riemannian manifolds, see e.g.
\cite{CaY}, \cite{BBBQ}, \cite{BG}, \cite{BW},\ \cite{Ga}, \cite{CKL}, \cite{CCF}, etc. Some of these works
concern Li-Yau type inequalities on pseudo-Hermitian manifolds. In
\cite{CKL}, Chang et al. derived a CR Li-Yau type estimate in terms of the
lower bound of pseudo-Hermitian Ricci curvature essentially for closed
Sasakian 3-manifolds. In \cite{BG}, Baudoin and Garofalo proved, among other
results, a CR Li-Yau type inequality on complete Sasakian manifolds under a
curvature dimension inequality. In \cite{CCL}, the authors announced a CR Li-Yau
gradient estimate by using a generalized curvature-dimension inequality and
the maximum principle in a closed pseudo-Hermitian manifold possibly with
nonvanishing torsion. Besides, they also established a Li-Yau type
inequality for the sum of squares of vector fields up to higher step on a
closed manifold, generalizing Cao-Yau's result (\cite{CaY}) for operators expressed as the
sum of squares of vector fields of step 2. However, we don't understand
their proof for the CR part. Anyhow Cao-Yau's
inequality in \cite{CaY} almost gave us in particular a Li-Yau type inequality for closed
pseudo-Hermitian manifolds (See Remark 3.1).
\par
Let's recall briefly Cao-Yau's work in \cite{CaY} as follows. Suppose $X_{1},....,$ $%
X_{n}$ are smooth vector fields on a closed manifold $M$ and 
\begin{equation}
L=\sum_{i=1}^{n}X_{i}^{2}-X_{0}   
\end{equation}%
with $X_{0}=\sum_{i=1}^{n}c_{i}X_{i}$, where $c_{i}$ are smooth functions on 
$M$. Suppose $X_{1},...,X_{n}$ satisfy the following conditions: for $1\leq
i,j,k\leq n$, $[X_{i},[X_{j},X_{k}]]$ can be expressed as linear
combinations of $X_{1},...,X_{n}$ and their brackets $
[X_{1},X_{2}],....,$\\$
[X_{n-1},X_{n}]$. Cao and Yau considered a positive
solution $u(x,t)$ of 
\begin{equation}
\left( L-\frac{\partial }{\partial t}\right) u=0  
\end{equation}%
on $M\times ( 0,\infty )$ and showed that there exists a constant $%
\delta _{0}>0$ such that for any $\delta >\delta _{0}$, $u$ satisfies 
\begin{equation}
\frac{1}{u^{2}}\sum_{i=1}^{n}|X_{i}u|^{2}-\delta \frac{X_{0}u}{u}-\delta 
\frac{u_{t}}{u}\leq \frac{C_{1}}{t}+C_{2} ,
\end{equation}%
where $C_{1}$ and $C_{2}$ are positive constants depending on $n$, $\delta
_{0}$, $\delta $ and $\{X_{i}\}$.
\par
This paper is devoted to establish a Li-Yau type inequality on a complete
pseudo-Hermitian manifold possibly with nonvanishing pseudo-Hermitian
torsion. The pseudo-Hermitian manifolds considered here are CR manifolds of
hypersurface type which admit positive definite pseudo-Hermitian structures
(see \S 2 for the detailed definition). Let $(M^{2m+1},HM,J,\theta )$ denote
a pseudo-Hermitian manifold of dimension $2m+1$. Here $(HM,J)$ is a CR
structure of type $(m,1)$, and $\theta $ is a pseudo-Hermitian structure on $%
M$. We find that a pseudo-Hermitian manifold carries a rich geometric
structure, including an almost complex structure $J$ on $HM$, the positive
definite Levi form $L_{\theta }$ on $HM$ induced from $\theta $ and $J$, the Webster metric (a Riemannian metric on $M$ extending $L_{\theta}$), the
Reeb vector field $\xi $ on $M$, the sub-Laplacian $\bigtriangleup _{b}$ (a
subelliptic differential operator) and the horizontal gradient operator $%
\nabla _{b}$ acting on functions. Note also that the pair $(HM,L_{\theta })$ is a $2$%
-step sub-Riemannian structure, which induces a Carnot--Carath\'{e}dory
distance $d_{cc}$ on $M$. These geometric data provide us a basis to
investigate Li-Yau type inequality on a pseudo-Hermitian manifold. We will
consider a positive solution of the following CR heat
equation 
\begin{equation}
\frac{\partial u}{\partial t}=\bigtriangleup _{b}u   
\end{equation}%
on a complete pseudo-Hermitian manifold, and establish a Li-Yau type
inequality for $u$. The main ingredients of Li-Yau's method \cite{LY} or
Cao-Yau's method (\cite{CaY}) involve the Bochner type formula, a parabolic
differential inequality for a suitable auxiliary function and the maximum
principle. For any smooth function $f$ on the pseudo-Hermitian manifold, one
has two CR Bochner formulas for $|\nabla _{b}f|^{2}$ and $f_{0}^{2}$
respectively, where $f_{0}=\xi (f)$ (see (2.8) and (2.9) in \S 2). Now set $%
f=\ln u$. Following Cao-Yau's idea, we will consider the auxiliary functions%
\begin{equation}
\mathcal{F}= t\left( |\nabla _{b}f|^{2}+t^{2\lambda -1}\left( 1+f_{0}^{2}\right)
^{\lambda }-\delta f_{t}\right) 
\end{equation}%
or 
\begin{equation}
\mathcal{G}=t\left( |\nabla _{b}f|^{2}+\left( 1+f_{0}^{2}\right) ^{\lambda }-\delta
f_{t}\right) 
\end{equation}%
according to the ranges of $t$. Some parabolic differential inequalities for $\mathcal{F}$ and $\mathcal{G}$ can be derived from the CR Bochner formulas. Following the technique in \cite{LY}, we may multiply $\mathcal{F}$ and $\mathcal{G}$ by a suitable cut-off function $\phi$
to localize the problem. By applying the maximum principle to $\phi \mathcal{F}$ and $\phi \mathcal{G}$,
and using the CR sub-Laplacian comparison theorem in \cite{CDRZ}, we are able to
establish the following local Li-Yau gradient estimate.
    
    ~\\
    $\mathbf{Theorem\ 1.1}$ Let $(M^{2m+1},HM,J,\theta)$ be a complete noncompact pseudo-Hermitian manifold with
    \begin{eqnarray}
        Ric_b+2(m-2)Tor_b \geq -k,\ and\  |A|,|\nabla_b A|\leq k_1,\nonumber
     \end{eqnarray}
     and $u$ be a positive solution of the CR heat equation
    \begin{eqnarray}
        \frac{\partial u}{\partial t}=\Delta_b u  \nonumber
     \end{eqnarray}
    on $B_p(2R)\times (0,\infty )\ with\ R\geq1 $, where $B_p(r)$ denotes the Riemannian ball of radius $r$ with respect to the Webster metric $g_{\theta}$. Then for any constant $\frac{1}{2}<\lambda<\frac{2}{3}$ and any constant $\delta>1+\frac{4}{m\lambda(2\lambda-1)}$, there exists a constant $C$ depending on $m,k,k_1,\lambda,\delta$, such that 
    \begin{eqnarray}
      \frac{|\nabla_b u|^2}{u^2}-\delta \frac{u_t}{u}\leq C(1+\frac{1}{t}+\frac{1}{R^\lambda}+\frac{1}{tR^\lambda})
    \end{eqnarray}
    on $B_p(R)\times (0,\infty )$.
   
   ~\par
    Letting $R\rightarrow \infty$ in Theorem 1.1, we get immediately the global Li-Yau type gradient estimate.
    
    ~\\
    $\mathbf{Theorem\ 1.2}$  Let $(M^{2m+1},HM,J,\theta)$ be a complete noncompact pseudo-Hermitian manifold with
    \begin{eqnarray}
        Ric_b+2(m-2)Tor_b\geq -k,\ and\   |A|,|\nabla_b A|\leq k_1,\nonumber
     \end{eqnarray}
     and $u$ be a positive solution of the heat equation
    \begin{eqnarray}
        \frac{\partial u}{\partial t}=\Delta_b u   \nonumber
     \end{eqnarray}
    on $M\times (0,\infty)$. Then for any constant $\frac{1}{2}<\lambda<\frac{2}{3}$ and any constant $\delta>1+\frac{4}{m\lambda(2\lambda-1)}$, there exists a constant $C$ depending on $m,k,k_1,\lambda,\delta$, such that 
    \begin{eqnarray}
      \frac{|\nabla_b u|^2}{u^2}-\delta \frac{u_t}{u}\leq C+\frac{C}{t}
    \end{eqnarray}
    on $M\times (0,\infty)$.
    
    ~\par
   As applications of the above gradient estimates, we give a Harnack inequality
for the positive solutions of the CR heat equation, and then obtain an upper bound
estimate for the heat kernel of the CR heat equation.
   
    ~\\
    $\mathbf{Theorem\ 1.3}$ Let $(M^{2m+1},HM,J,\theta)$ be a complete noncompact pseudo-Hermitian manifold with
    \begin{eqnarray}
       Ric_b+2(m-2)Tor_b \geq -k,\ and\  |A|,|\nabla_b A|\leq k_1,\nonumber
     \end{eqnarray}
     and $u$ be a positive solution of the heat equation
    \begin{eqnarray}
        \frac{\partial u}{\partial t}=\Delta_b u  \nonumber
     \end{eqnarray}
    on $M\times (0,\infty)$. Then for any constant $\frac{1}{2}<\lambda<\frac{2}{3}$ and any constant $\delta>1+\frac{4}{m\lambda(2\lambda-1)}$, there exists a constant $C$ which is given by Theorem 1.2 such that for any $0<t_1<t_2$ and $x,y\in M$, we have
    \begin{eqnarray}
      u(x,t_1)\leq u(y, t_2)(\frac{t_2}{t_1})^{\frac{C}{\delta}}exp(\frac{C}{\delta}(t_2-t_1)+\frac{\delta d_{cc}^2(x,y)}{4(t_2-t_1)} ).
    \end{eqnarray} 
   
   ~\\
    $\mathbf{Theorem\ 1.4}$ Let $(M^{2m+1},HM,J,\theta)$ be a complete noncompact pseudo-Hermitian manifold with
    \begin{eqnarray}
        Ric_b+2(m-2)Tor_b \geq -k,\ and\  |A|,|\nabla_b A|\leq k_1,\nonumber
     \end{eqnarray}
    and $H(x,y,t)$ be the heat kernel of (1.4). Then for any constants $\frac{1}{2}<\lambda<\frac{2}{3}$, $\delta>1+\frac{4}{m\lambda(2\lambda-1)}$ and  $0<\epsilon <1$, there exists constants $C^{'}\ and\ C^{''}$ depending on $m,k,k_1,\lambda,\delta,\epsilon$, such that $H(x,y,t)$ satisfies 
    \begin{align}
          H(x,y,t)\leq C^{'} [Vol(B_{cc}(x,\sqrt{t}))]^{-\frac{1}{2}} [Vol(B_{cc}(y,\sqrt{t}))]^{-\frac{1}{2}} exp(C^{''}\epsilon t-\frac{d^2_{cc}(x,y)}{(4+\epsilon)t}),
    \end{align}
    where the $B_{cc}(x,r) $ is the ball with respect to Carnot-Carath$\acute{e}$odory distance. The constant $C^{'}\rightarrow \infty$ as $\epsilon \rightarrow 0$.




\section{CR Bochner formulas on pseudo-Hermitian manifolds}
\label{2}

In this section we introduce some basic notations in pseudo-Hermitian geometry (cf. \cite{DT,We,Ta} for details), and then give the CR Bochner formulas for functions on a pseudo-Hermitian manifold. Next, we will derive parabolic
differential inequalities for the auxiliary functions $\mathcal{F}$ and $\mathcal{G}$.
\par 
 Let $M^{2m+1}$ be a real $2m + 1$ dimensional orientable $C^\infty$ manifold. A CR structure on $M$ is a complex subbundle $H^{1,0}M $ of $TM\otimes \mathbb{C}$ satisfying 
\begin{equation}
H^{1,0}M\cap H^{0,1}M=\{0\} ,\ \  [\Gamma(H^{1,0}M),\Gamma(H^{1,0}M)]\subseteq \Gamma(H^{1,0}M)   
\end{equation}
where $H^{0,1}M=\overline{H^{1,0}M}$. Equivalently, the CR structure may also be described by the real bundle $HM=Re\{H^{1,0}M\oplus H^{0,1}M\}$ and an almost complex structure $J$ on $HM$, where $J(X+\overline{X})=\sqrt{-1}(X-\overline{X})$ for any $X\in H^{1,0}M$. Then $(M,HM,J)$ is said to be a CR manifold.
\par 
We denote by $E$ the conormal bundle of $HM$ in $T^*M$, whose fiber at each
point $x\in M$ is given by
\begin{eqnarray}
E_x=\{\omega\in T_x^*M|\omega(H_xM)=0\}.
\end{eqnarray}
It turns out that $E$ is a trivial line bundle. Therefore there exist globally defined nowhere vanishing sections $\theta \in \Gamma(E)$.  A section $\theta \in \Gamma(E\backslash \{0\})$ is called a pseudo-Hermitian structure on $M$. The Levi form $L_\theta$ of a pseudo-Hermitian structure $\theta$ is defined by
\begin{eqnarray}
L_\theta(X,Y)=d\theta(X,JY) \nonumber
\end{eqnarray}
for any $X,Y\in HM$. The integrability condition in (2.1) implies that $L_\theta $ is $J$-invariant, and thus symmetric. When $L_\theta$ is positive definite on $HM$ for some $\theta$,
then $(M, HM, J)$ is said to be strictly pseudoconvex. From now on, we will always assume that $(M, HM, J)$ is a strictly pseudoconvex CR manifold endowed
with $\theta$, such that $L_\theta$ is positive definite. Then the quadruple $(M, HM, J, \theta)$ is referred to as a pseudo-Hermitian manifold.
\par
For a pseudo-Hermitian manifold $(M, HM, J, \theta)$, due to the positivity of $L_\theta$, we have a sub-Riemannian structure $(HM,L_\theta)$ of step-2 on $M$. We say that a Lipschitz curve $\gamma:[0,l]\rightarrow M$ is horizontal if $\gamma^{'} \in H_{\gamma(t)}M$ a.e. in $[0,l]$. For any two points $p,q\in M$, by the well-known theorem of Chow-Rashevsky(\cite{Ch,Ra} ),  there always exist such horizontal curves joining $p$ and $q$. Therefore we
may define the Carnot-Carath$\acute{e}$odory distance as follows:
\begin{eqnarray}
d_{cc}(p,q)= inf\{ \int_0^l \sqrt{L_\theta(\gamma^{'},\gamma^{'})} dt \ | \  \gamma \in \Gamma(p,q)\},
\nonumber
\end{eqnarray}
where $\Gamma(p,q) $ denotes the set of all horizontal curves joining $p$ and $q$. Clearly $d_{cc}$ induces to a metric space structure on $M$, in which its metric ball centered at $x$ with radius $r$ is given by
\begin{eqnarray}
B_{cc}(x,r)=\{y\in M \ | \  d_{cc}(y,x)<r\}. \nonumber
\end{eqnarray}
\par
For a pseudo-Hermitian manifold $(M, HM, J, \theta)$, it is clear that $\theta$ is a contact form on $M$. Consequently there exists a unique vector field $\xi$ such that
\begin{eqnarray}
\theta(\xi)=1,\ d\theta(\xi,\cdot)=0. 
\end{eqnarray}
This vector field $\xi$ is called the Reeb vector field. From (2.2) and (2.3), it is easy to see that $TM$ admits the following direct sum decomposition
\begin{eqnarray}
TM=HM \oplus R\xi ,
\end{eqnarray}
which induces a natural projection $\pi_b:TM\rightarrow HM$. In terms of $\theta$ and the decomposition (2.4), the Levi form $L_\theta$ can be extended to a Riemannian metric
\begin{eqnarray}
g_\theta=L_\theta+\theta\otimes \theta , \nonumber
\end{eqnarray}
which is called the Webster metric. We will denote by $r$ the corresponding Riemannian distance and by $B_p(R)$ the Riemannian ball of radius $R$ centered at $p$. One may extend the complex structure $J$
on $HM$ to an endomorphism of $TM$, still denoted by $J$, by requiring
\begin{eqnarray}
J\xi=0. \nonumber
\end{eqnarray}
It is known that there exists a canonical connection $\nabla$ on a pseudo-Hermitian manifold, called the Tanaka-Webster connection (cf. \cite{DT,Ta,We}), such that 
\begin{eqnarray}
&1.&\ \nabla_X\Gamma(HM)\subseteq \Gamma(HM),\ for\ any\ X\in \Gamma(TM); \nonumber\\
&2.&\ \nabla g_\theta=0\ and\  \nabla J=0; \nonumber\\
&3.&\ T_\nabla(X,Y)=2d\theta(X,Y)\xi\ and\ T_\nabla(\xi,JX)+JT_\nabla(\xi,X)=0, \nonumber\\
&&for \ any\  X,Y\in HM,\ where\ T_\nabla \ denotes\ the\ torsion\ of
\nonumber\\
&&the\ connection\ \nabla. \nonumber
\end{eqnarray}
The pseudo-Hermitian torsion of $\nabla$ is an important pseudo-Hermitian invariant, which is an $HM$-valued 1-form defined by
\begin{eqnarray}
\tau(X)=T_\nabla (\xi,X) \nonumber
\end{eqnarray}
for any $X\in TM$. Note that is $\tau$ trace-free and self-adjoint with respect to the
Webster metric $g_\theta$ (cf. \cite{DT}). Set $A(X,Y)=g_\theta(T_\nabla (\xi,X),Y )$ for any $X,Y\in TM$, then we have
\begin{eqnarray}
A(X,Y)=A(Y,X) .
\end{eqnarray}
We say that $M$ is Sasakian if $\tau=0$ (or  equivalently, $A=0$).

\par 
Let $(M,HM,J,\theta)$ be a complete pseudo-Hermitian manifold of dimension $2m+1$. We choose a local orthonormal frame field $\{ e_A\}_{A=0}^{2m}=\{\xi,e_1,\cdots,e_m,\\
e_{m+1},\cdots,e_{2m} \} $ with respect to the Webster metric $g_\theta$ such that
\begin{eqnarray}
\{e_{m+1},\cdots,e_{2m}\}=\{ Je_1,\cdots,Je_m\} .\nonumber
\end{eqnarray}
Set
\begin{eqnarray}
\eta_\alpha=\frac{1}{\sqrt{2}}(e_\alpha-\sqrt{-1}J e_\alpha),\quad  \eta_{\bar{\alpha}}=\frac{1}{\sqrt{2}}(e_\alpha+\sqrt{-1}J e_\alpha),\ (\alpha=1,\cdots,m). \nonumber
\end{eqnarray}
Then $\{\eta_\alpha\}_{\alpha=1}^m$ is a unitary frame field of $H^{1,0}M$ with respect to $g_\theta$. Let $\{ \theta^1,\cdots,\theta^m\}$ be the dual frame field of $\{\eta_\alpha \}_{\alpha=1}^m$. According to the property 3 of the Tanaka-Webster connection, one may write
\begin{eqnarray}
\tau&=&\tau^\alpha\eta_\alpha+\tau^{\bar{\alpha}}\eta_{\bar{\alpha}} 
\nonumber\\
&=&A_{\bar{\beta}}^\alpha \theta^{\bar{\beta}}\otimes\eta_\alpha+A_\beta^{\bar{\alpha}}\theta^\beta\otimes\eta_{\bar{\alpha}} .\nonumber
\end{eqnarray}
We will also write $A_{\alpha\beta}=A_{\beta}^{\bar{\alpha}}$ and $A_{\bar{\alpha}\bar{\beta}}=A_{\bar{\beta}}^\alpha$. Then (2.5) means that $A_{\alpha\beta}=A_{\beta\alpha}$ and $A_{\bar{\alpha}\bar{\beta}}=A_{\bar{\beta}\bar{\alpha}}$. From \cite{We},  we have the following structure equations of the Tanaka-Webster connection $\nabla$:
\begin{eqnarray}
d\theta&=&2\sqrt{-1}\theta^\alpha \wedge \theta^{\bar{\alpha}}, \nonumber\\
d\theta^{\alpha}&=&\theta^\beta\wedge \theta^\alpha_\beta+A_{\bar{\alpha}\bar{\beta}}\theta\wedge\theta^\beta,\\
d\theta^\alpha_\beta&=&\theta^\gamma_\beta\wedge\theta^\alpha_\gamma+\Pi^\alpha_\beta \nonumber
\end{eqnarray}
with
\begin{eqnarray}
\Pi^\alpha_\beta=2\sqrt{-1}(\theta^\alpha\wedge\tau^{\bar{\beta}}-\tau^\alpha \wedge \theta^{\bar{\beta}})+R^\alpha_{\beta\lambda\bar{\mu}}\theta^\lambda\wedge\theta^{\bar{\mu}}+W^\alpha_{\beta\bar{\gamma}}\theta\wedge\theta^{\bar{\gamma}}-W^\alpha_{\beta\gamma}\theta\wedge\theta^{{\gamma}}, \nonumber
\end{eqnarray}   
where $W^\alpha_{\beta\bar{\gamma}}=A^\alpha_{\bar{\gamma}, \beta},\ W^\alpha_{\beta\gamma}=A^{\bar{\gamma}}_{\beta,\bar{\alpha}}$ are the are the covariant derivatives of $A$, and $R_{\beta\lambda\bar{\mu}}^\alpha$ are the components of curvature tensor of the Tanaka-Webster connection. Set
\begin{eqnarray}
R_{\alpha\bar{\beta}}=R_{\gamma\alpha\bar{\beta}}^\gamma ,\nonumber
\end{eqnarray}
then $R_{\alpha\bar{\beta}}=R_{\bar{\beta}\alpha}$ (cf. \cite{DT}). For any $X=a^\alpha \eta_\alpha+b^{\bar{\alpha}}\eta_{\bar{\alpha}}$ and $Y=c^\beta\eta_\beta+d^{\bar{\beta}}\eta_{\bar{\beta}} \in HM\otimes \mathbb{C}$, we define
\begin{eqnarray}
  Ric_b(X,Y)=R_{\alpha\bar{\beta}}a^\alpha d^{\bar{\beta}}+R_{\bar{\alpha}\beta}b^{\bar{\alpha}}c^\beta ,\nonumber
\end{eqnarray}  
whose components are given by
\begin{eqnarray}
  Ric_b(\eta_\alpha,\eta_{\bar{\beta}})&=&R_{\alpha\bar{\beta}},\  Ric_b(\eta_{\bar{\alpha}},\eta_{\beta})=R_{\bar{\alpha}\beta} ,\nonumber\\
  Ric_b(\eta_\alpha,\eta_{{\beta}})&=& Ric_b(\eta_{\bar{\alpha}},\eta_{\bar{\beta}})=0. \nonumber
\end{eqnarray} 
The 2-tensor $Ric_b$ will be referred to as the pseudo-Hermitian Ricci tensor. For any $X=X^\alpha \eta_\alpha+X^{\bar{\alpha}}\eta_{\bar{\alpha}}$ and $Y=Y^\beta\eta_\beta+Y^{\bar{\beta}}\eta_{\bar{\beta}} \in HM\otimes \mathbb{C}$, we introduce
\begin{eqnarray}
 Tor_b(X,Y)&=& A(X,JY) \nonumber\\
 &=&\sqrt{-1}A(X^\alpha\eta_\alpha+X^{\bar{\alpha}}\eta_{\bar{\alpha}},Y^\beta\eta_\beta-Y^{\bar{\beta}}\eta_{\bar{\beta}} ) \nonumber\\
 &=& \sqrt{-1}( A_{\alpha\beta}X^\alpha Y^\beta -A_{\bar{\alpha}\bar{\beta}}X^{\bar{\alpha}}Y^{\bar{\beta}}). \nonumber
\end{eqnarray}  
Clearly both $Ric_b$ and $Tor_b$ are real symmetric, fiberwise 2-tensors on $HM$.
\par
For a $C^2$ function $f:M\rightarrow R$, its differential $df$ and gradient $\nabla f$ can be expressed as
\begin{eqnarray}
df=f_0\theta+f_\alpha\theta^\alpha+f_{\bar{\alpha}}\theta^{\bar{\alpha}} \nonumber
\end{eqnarray}
and
\begin{eqnarray}
\nabla f=f_0\xi+ f_{\bar{\alpha}}\eta_\alpha+f_\alpha \eta_{\bar{\alpha}} ,\nonumber
\end{eqnarray}
where $f_0=\xi(f),f_\alpha=\eta_\alpha(f),f_{\bar{\alpha}}=\eta_{\bar{\alpha}}(f)$. Then the horizontal gradient of $f$ is given by
\begin{eqnarray}
\nabla_b f=f_{\bar{\alpha}}\eta_\alpha+f_\alpha \eta_{\bar{\alpha}} .\nonumber
\end{eqnarray}
Let $\nabla d f$ be the covariant derivative of the differential $df\in \Gamma(T^*M)$ with respect to the Tanaka-Webster connection. Then $\nabla d f$ may be expressed as 
\begin{eqnarray}
 \nabla d f&=& f_{\alpha\beta}\theta^\alpha\otimes\theta^\beta +f_{\alpha\bar{\beta}}\theta^\alpha \otimes \theta^{\bar{\beta}}+ f_{\bar{\alpha}\beta}\theta^{\bar{\alpha}}\otimes \theta^\beta+f_{\bar{\alpha}\bar{\beta}}\theta^{\bar{\alpha}}\otimes \theta^{\bar{\beta}} 
 \nonumber\\
 &&+f_{0\alpha}\theta\otimes\theta^\alpha+f_{0\bar{\alpha}}\theta\otimes\theta^{\bar{\alpha}}+f_{\alpha0}\theta^\alpha\otimes\theta+f_{\bar{\alpha}0}\theta^{\bar{\alpha}}\otimes \theta .\nonumber
\end{eqnarray}
The following communication relations are known  (see, e.g., Chapter 9 in \cite{DT}, or \S 3 in \cite{CDRY}):
 \begin{eqnarray}
      f_{\alpha\beta}=f_{\beta\alpha},\quad  f_{\alpha\bar{\beta}}-f_{\bar{\beta}\alpha}=2\sqrt{-1}f_0\delta_\alpha^\beta,\quad f_{0\alpha}-f_{\alpha 0}=f_{\bar{\beta}}A_{\alpha}^{\bar{\beta}}.
 \end{eqnarray}
 The horizontal Hessian of $f$ is defined by
\begin{eqnarray}
     Hess_b(f)&=& (\nabla d f)(\pi_b,\pi_b) \nonumber\\
     &=& f_{\alpha\beta}\theta^\alpha\otimes\theta^\beta +f_{\alpha\bar{\beta}}\theta^\alpha \otimes \theta^{\bar{\beta}}+ f_{\bar{\alpha}\beta}\theta^{\bar{\alpha}}\otimes \theta^\beta+f_{\bar{\alpha}\bar{\beta}}\theta^{\bar{\alpha}}\otimes \theta^{\bar{\beta}}.\nonumber
 \end{eqnarray}
 Consequently
    \begin{eqnarray}
      |\nabla_b f|^2=2f_\alpha f_{\bar{\alpha}}, \quad  |Hess_b(f)|^2=2(f_{\alpha\beta}f_{\bar{\alpha}\bar{\beta}}+f_{\alpha\bar{\beta}}f_{\bar{\alpha}\beta}) .\nonumber
    \end{eqnarray}
 The sub-Laplacian of $f$ is defined by
 \begin{eqnarray}
      \Delta_b f&=& tr\{Hess_b(f)\} \nonumber\\
      &=&f_{\alpha\bar{\alpha}}+f_{\bar{\alpha}\alpha}. \nonumber
 \end{eqnarray}
 From \cite{GL}, \cite{Gr} (see also \cite{DT} and \cite{CDRY}), we have the following Bochner formulas
   \begin{eqnarray}
     \frac{1}{2}\Delta_b |\nabla_b f|^2&=& 2(f_{\alpha\beta}f_{\bar{\alpha}\bar{\beta}}+f_{\alpha\bar{\beta}}f_{\bar{\alpha}\beta})+f_{\bar{\alpha}}(f_{\beta\bar{\beta}}+f_{\bar{\beta}\beta})_{\alpha}+f_\alpha (f_{\beta\bar{\beta}}+f_{\bar{\beta}\beta})_{\bar{\alpha}} \nonumber\\
     &&+2R_{\alpha\bar{\beta}}f_{\bar{\alpha}}f_\beta+2\sqrt{-1}(m-2)(A_{\alpha\beta}f_{\bar{\alpha}}f_{\bar{\beta}}-A_{\bar{\alpha}\bar{\beta}}f_\alpha f_\beta ) \nonumber\\
     &&+4\sqrt{-1}(f_{\bar{\alpha}}f_{ 0\alpha}-f_{\alpha}f_{0\bar{\alpha}}) 
     \\
     &=&|Hess_b(f)|^2+\langle \nabla_b f,\nabla_b \Delta_b f\rangle +4\langle J\nabla_b f,\nabla_b f_0\rangle
     \nonumber\\
     &&+(Ric_b+2(m-2)Tor_b)(\nabla_b f,\nabla_b f) \nonumber
    \end{eqnarray}
    and 
    \begin{eqnarray}
      \frac{1}{2}\Delta_b f_0^2&=&|\nabla_b f_0|^2+f_0(\Delta_b f)_0
      \nonumber\\
      &&+2f_0(f_\beta A_{\bar{\beta}\bar{\alpha},\alpha}+f_{\bar{\beta}}A_{\beta\alpha,\bar{\alpha}}+f_{\beta\alpha}A_{\bar{\beta}\bar{\alpha}}+f_{\bar{\beta}\bar{\alpha}}A_{\beta\alpha}) 
      \\
      &=& |\nabla_b f_0|^2+f_0(\Delta_b f)_0+2f_0 \mathrm{Im} Qf, \nonumber
    \end{eqnarray}
   where $Q$ is the purely holomorphic second-order operator defined by (\cite{GL})
    \begin{eqnarray}
      Qf=2\sqrt{-1}(A_{\bar{\beta}\bar{\alpha}}f_\beta)_\alpha .\nonumber
    \end{eqnarray}
    Note that the coefficient before the 'mixed term' $\langle J\nabla_b f,\nabla_b f_0\rangle$ in (2.8) is slight different from that in \cite{Gr}.
    
    ~\\
    $\mathbf{Lemma\ 2.1}$  Let $(M^{2m+1},HM,J,\theta)$ be a pseudo-Hermitian manifold and $u$ be a positive solution of the CR heat equation (1.4). Set $f=\mathrm{ln}\ u$. Then for any $0<\lambda\leq 1$,
    we have 
    \begin{eqnarray}
        (\Delta_b-\partial_t)|\nabla_b f|^2&\geq & \frac{1}{m}(\Delta_b f)^2+4mf^2_0+4f_{\alpha\beta}f_{\bar{\alpha}\bar{\beta}}
                \nonumber\\
        &&-2\langle \nabla_b |\nabla_b f|^2,\nabla_b f \rangle +8 \langle \nabla_b f_0, J\nabla_b f \rangle
        \nonumber\\
        &&+2(Ric_b +2(m-2)Tor_b)(\nabla_b f,\nabla_b f) 
     \end{eqnarray}
     and 
     \begin{eqnarray}
        (\Delta_b-\partial_t)(1+f^2_0)^\lambda
        &\geq&2\lambda(2\lambda-1) (1+f^2_0)^{\lambda-1}|\nabla_b f_0|^2
        \nonumber\\
        &&-2\lambda (1+f^2_0)^{\lambda-1}(\langle \nabla_b f^2_0, \nabla_b f \rangle -2f_0A(\nabla_bf, \nabla_b f)) 
        \nonumber\\
        &&+4\lambda (1+f^2_0)^{\lambda-1}f_0\mathrm{Im} Qf .
     \end{eqnarray}

     $\mathbf{Proof}$ A direct computation gives
     \begin{eqnarray}
     (\Delta_b-\partial_t) f=-|\nabla_b f|^2 .
    \end{eqnarray} 
    Using (2.7), we may estimate the term in $|Hess_b(f)|^2$:
    \begin{eqnarray}
         f_{\alpha\bar{\beta}}f_{\bar{\alpha}\beta}
        &\geq& \sum\limits_{\alpha=1}^m f_{\alpha\bar{\alpha}}f_{\bar{\alpha}\alpha}=\frac{1}{4}\sum\limits_{\alpha=1}^m (|f_{\alpha\bar{\alpha}}-f_{\bar{\alpha}\alpha}|^2+|f_{\alpha\bar{\alpha}}+f_{\bar{\alpha}\alpha}|^2   ) \nonumber\\
        &\geq& \frac{1}{4m}|\sum\limits_{\alpha=1}^m( f_{\alpha\bar{\alpha}}+f_{\bar{\alpha}\alpha}) |^2+\frac{1}{4}\sum\limits_{\alpha=1}^m|f_{\alpha\bar{\alpha}}-f_{\bar{\alpha}\alpha}|^2 \nonumber\\
        &=&\frac{1}{4m}(\Delta_b f)^2+mf_0^2.
     \end{eqnarray}
     Then (2.10) follows immediately from (2.8), (2.12) and (2.13). From (2.9), we get 
    \begin{eqnarray}
      (\Delta_b-\partial_t) f_0^2
     =  2|\nabla_b f_0|^2+2f_0(\Delta_b f -\partial_t f)_0+4f_0 \mathrm{Im} Qf .
    \end{eqnarray}
    Using (2.14), we derive that
     \begin{eqnarray}
        &&(\Delta_b-\partial_t)(1+f^2_0)^\lambda \nonumber\\
        &=&2\lambda(\lambda-1)(1+f^2_0)^{\lambda-2}(f^2_0)_{\alpha}(f^2_0)_{\Bar{\alpha}}+\lambda (1+f^2_0)^{\lambda-1}(f^2_0)_{\alpha\Bar{\alpha}+\Bar{\alpha}\alpha} \nonumber\\
        &&-\partial_t(1+f^2_0)^\lambda 
        \nonumber\\
        &=&4\lambda (\lambda-1)(1+f^2_0)^{\lambda-2}f^2_0|\nabla_b f_0|^2+\lambda (1+f^2_0)^{\lambda-1}(\Delta_b-\partial_t)(f^2_0) 
        \nonumber\\
        &=&4\lambda (\lambda-1)(1+f^2_0)^{\lambda-2}f^2_0|\nabla_b f_0|^2+2\lambda (1+f^2_0)^{\lambda-1}|\nabla_b f_0|^2 
        \nonumber\\
        &&+2\lambda (1+f^2_0)^{\lambda-1}f_0\nabla_\xi(\Delta_b f -f_t)
        \nonumber\\
        &&+4\lambda (1+f^2_0)^{\lambda-1}f_0\mathrm{Im} Qf
        \nonumber\\
        &=&2\lambda (1+f^2_0)^{\lambda-2}|\nabla_b f_0|^2((2\lambda-1)f^2_0+1) 
        \nonumber\\
        &&-2\lambda (1+f^2_0)^{\lambda-1}(\langle \nabla_b f^2_0, \nabla_b f \rangle -2f_0A(\nabla_bf, \nabla_b f)) 
        \nonumber\\
        &&+4\lambda (1+f^2_0)^{\lambda-1}f_0\mathrm{Im} Qf
        \nonumber\\
        &\geq&2\lambda(2\lambda-1) (1+f^2_0)^{\lambda-1}|\nabla_b f_0|^2
        \nonumber\\
        &&-2\lambda (1+f^2_0)^{\lambda-1}(\langle \nabla_b f^2_0, \nabla_b f \rangle -2f_0A(\nabla_bf, \nabla_b f)) 
        \nonumber\\
        &&+4\lambda (1+f^2_0)^{\lambda-1}f_0\mathrm{Im} Qf .
        \nonumber
     \end{eqnarray}
    This completes the proof of Lemma 2.1. 
    {\qed}
    
\section{Li-Yau type gradient estimates}
\label{3}
\par
In this section, we derive the CR version of Li-Yau type gradient estimates. Inspired by Cao-Yau's method in \cite{CaY} for weakly elliptic operators on closed manifolds, we consider the auxiliary functions $\mathcal{F}$ and $\mathcal{G}$ on complete pseudo-Hermitian manifolds, trying to establish the gradient estimates according to two cases: $0< t\leq 1$ and $t\geq1$. \par 
Let $u$ be a positive solution of (1.4). Recall that, if $f=\mathrm{ln}\ u$, then
\begin{eqnarray}
(\Delta_b-\frac{\partial}{\partial t})f=-|\nabla_b f|^2.
\end{eqnarray}
As mentioned in Introduction, we consider the following auxiliary functions
    \begin{equation}
 \mathcal{F}=t\left( |\nabla _{b}f|^{2}+t^{2\lambda -1}\left( 1+f_{0}^{2}\right)
^{\lambda }-\delta f_{t}\right)=\mathcal{F}_1+\mathcal{F}_2  
\end{equation}%
and 
\begin{equation}
\mathcal{G}=t\left( |\nabla _{b}f|^{2}+\left( 1+f_{0}^{2}\right) ^{\lambda }-\delta
f_{t}\right) ,
\end{equation}%
     where
    \begin{eqnarray}
       \mathcal{F}_1=t(|\nabla_b f|^2-\delta f_t),\  \mathcal{F}_2=t^{2\lambda}(1+f_0^2)^\lambda ,\nonumber
     \end{eqnarray}
     and $\delta>1$ is a constant. Using Lemma 2.1, we derive the following inequalities
     \begin{eqnarray}
        (\Delta_b-\partial_t)\mathcal{F}&=&(\Delta_b-\partial_t)\mathcal{F}_1+(\Delta_b-\partial_t)\mathcal{F}_2
        \nonumber\\
        &=&t(\Delta_b-\partial_t)(|\nabla_b f|^2-\delta f_t)-(|\nabla_b f|^2-\delta f_t)
        \nonumber\\
        &&+t^{2\lambda}(\Delta_b-\partial_t)(1+f_0^2)^\lambda-2\lambda t^{2\lambda-1} (1+f_0^2)^\lambda
        \nonumber\\
        &\geq &
        -\frac{\mathcal{F}_1}{t}-\frac{2\lambda \mathcal{F}_2}{t}+t\{ \frac{1}{m}(\Delta_b f)^2+4mf^2_0+4f_{\alpha\beta}f_{\Bar{\alpha}\Bar{\beta}}
        \nonumber\\
        &&+8 \langle \nabla_b f_0, J\nabla_b f \rangle +2(Ric_b +2(m-2)Tor_b)(\nabla_b f,\nabla_b f) \}
        \nonumber\\
       &&+t^{2\lambda} \{2\lambda(2\lambda-1) (1+f^2_0)^{\lambda-1}|\nabla_b f_0|^2
        \nonumber\\
        &&+4\lambda (1+f^2_0)^{\lambda-1}f_0A(\nabla_bf, \nabla_b f)
        \nonumber\\
        &&+4\lambda (1+f^2_0)^{\lambda-1}f_0\mathrm{Im} Qf\}
        \nonumber\\
        &&-2t\langle \nabla_b f,\nabla_b |\nabla_b f|^2\rangle +2\delta t\langle \nabla_b f,\nabla_b f_t\rangle
        \nonumber\\
        &&-2\lambda t^{2\lambda}(1+f_0^2)^{\lambda-1} \langle \nabla_b f,\nabla_b f_0^2\rangle
        \nonumber\\
        &\geq &
        -\frac{\mathcal{F}_1}{t}-\frac{2\lambda \mathcal{F}_2}{t}+t\{ \frac{1}{m}(\Delta_b f)^2+4mf^2_0+4f_{\alpha\beta}f_{\Bar{\alpha}\Bar{\beta}}
        \nonumber\\
        &&+8 \langle \nabla_b f_0, J\nabla_b f \rangle +2(Ric_b +2(m-2)Tor_b)(\nabla_b f,\nabla_b f) \}
        \nonumber\\
       &&+t^{2\lambda} \{2\lambda(2\lambda-1) (1+f^2_0)^{\lambda-1}|\nabla_b f_0|^2
        \nonumber\\
        &&+4\lambda (1+f^2_0)^{\lambda-1}f_0A(\nabla_bf, \nabla_b f)
        \nonumber\\
        &&+4\lambda (1+f^2_0)^{\lambda-1}f_0 \mathrm{Im} Qf\}-2\langle \nabla_b f,\nabla_b \mathcal{F}\rangle 
      \end{eqnarray}
      and 
       \begin{eqnarray}
        (\Delta_b-\partial_t)\mathcal{G}&=&
        t(\Delta_b-\partial_t)(|\nabla_b f|^2+( 1+f_{0}^{2}) ^{\lambda }-\delta f_t)
        \nonumber\\
        &&-(|\nabla_b f|^2+( 1+f_{0}^{2}) ^{\lambda }-\delta f_t)
        \nonumber\\
        &\geq &
        -\frac{\mathcal{G}}{t}+t\{ \frac{1}{m}(\Delta_b f)^2+4mf^2_0+4f_{\alpha\beta}f_{\Bar{\alpha}\Bar{\beta}}
        \nonumber\\
        &&+8 \langle \nabla_b f_0, J\nabla_b f \rangle +2(Ric_b +2(m-2)Tor_b)(\nabla_b f,\nabla_b f) 
        \nonumber\\
       &&+2\lambda(2\lambda-1) (1+f^2_0)^{\lambda-1}|\nabla_b f_0|^2
        \nonumber\\
        &&+4\lambda (1+f^2_0)^{\lambda-1}f_0A(\nabla_bf, \nabla_b f)
        \nonumber\\
        &&+4\lambda (1+f^2_0)^{\lambda-1}f_0\mathrm{Im} Qf\}
        \nonumber\\
        &&-2t\langle \nabla_b f,\nabla_b |\nabla_b f|^2\rangle +2\delta t\langle \nabla_b f,\nabla_b f_t\rangle
        \nonumber\\
        &&-2\lambda t(1+f_0^2)^{\lambda-1} \langle \nabla_b f,\nabla_b f_0^2\rangle
        \nonumber\\
        &\geq &
        -\frac{\mathcal{G}}{t}+t\{ \frac{1}{m}(\Delta_b f)^2+4mf^2_0+4f_{\alpha\beta}f_{\Bar{\alpha}\Bar{\beta}}
        \nonumber\\
        &&+8 \langle \nabla_b f_0, J\nabla_b f \rangle +2(Ric_b +2(m-2)Tor_b)(\nabla_b f,\nabla_b f) 
        \nonumber\\
       &&+2\lambda(2\lambda-1) (1+f^2_0)^{\lambda-1}|\nabla_b f_0|^2
        \nonumber\\
        &&+4\lambda (1+f^2_0)^{\lambda-1}f_0A(\nabla_bf, \nabla_b f)
        \nonumber\\
        &&+4\lambda (1+f^2_0)^{\lambda-1}f_0 \mathrm{Im} Qf\}-2\langle \nabla_b f,\nabla_b \mathcal{G}\rangle .
      \end{eqnarray}
    Hence we have the following Lemma 3.1.

    ~\\
     $\mathbf{Lemma\ 3.1}$ Let $(M^{2m+1},\theta)$ be a complete noncompact pseudo-Hermitian manifold with
     \begin{eqnarray}
       Ric_b +2(m-2)Tor_b\geq -k\ \  and\ \  |A|,|\nabla_b A|\leq k_1,\nonumber
     \end{eqnarray}
     then for $\frac{1}{2}<\lambda<1$, we have 
      \begin{eqnarray}
        (\Delta_b-\partial_t)\mathcal{F}&\geq &-\frac{ \mathcal{F}_1}{t}-\frac{2\lambda \mathcal{F}_2}{t}-2\langle \nabla_b f,\nabla_b \mathcal{F} \rangle +t\{ \frac{1}{m}(\Delta_b f)^2+4m f^2_0 
        \nonumber\\
        &&-\frac{8}{\lambda(2\lambda-1)}(1+f_0^2)^{1-\lambda}|\nabla_bf |^2\cdot t^{1-2\lambda} -2k|\nabla_b f|^2
        \nonumber\\
        &&-2k_1|\nabla_b f|^2-4k_1\lambda (1+f^2_0)^{\lambda-1}|f_0|\cdot|\nabla_b f|^2 \cdot t^{2\lambda-1}
        \nonumber\\
        &&-2(k_1+k_1^2)\lambda^2(1+f_0^2)^{2\lambda-2}|f_0^2|\cdot t^{4\lambda-2} \}
     \end{eqnarray}
     for $0<t\leq 1 $ and 
     \begin{eqnarray}
        (\Delta_b-\partial_t)\mathcal{G}&\geq &-\frac{\mathcal{G}}{t}-2\langle \nabla_b f,\nabla_b \mathcal{G} \rangle +t\{ \frac{1}{m}(\Delta_b f)^2+4m f^2_0 
        \nonumber\\
        &&-\frac{8}{\lambda(2\lambda-1)}(1+f_0^2)^{1-\lambda}|\nabla_bf |^2 -2k|\nabla_b f|^2
        \nonumber\\
        &&-2k_1|\nabla_b f|^2-4k_1\lambda (1+f^2_0)^{\lambda-1}|f_0|\cdot|\nabla_b f|^2 
        \nonumber\\
        &&-2(k_1+k_1^2)\lambda^2(1+f_0^2)^{2\lambda-2}|f_0^2|\}
     \end{eqnarray}
     for $t\geq 1$.
     ~\\
     $\mathbf{Proof}$ From (3.4) and the condition, we find that
    \begin{eqnarray}
        (\Delta_b-\partial_t)\mathcal{F}&\geq &-\frac{\mathcal{F}_1}{t}-\frac{2\lambda \mathcal{F}_2}{t}+t\{ \frac{1}{m}(\Delta_b f)^2+4m f^2_0+4f_{\alpha\beta}f_{\Bar{\alpha}\Bar{\beta}}
        \nonumber\\
        &&+8 \langle \nabla_b f_0, J\nabla_b f \rangle -2k|\nabla_b f|^2
        \nonumber\\
       &&+2\lambda(2\lambda-1) (1+f^2_0)^{\lambda-1}|\nabla_b f_0|^2\cdot t^{2\lambda-1} 
        \nonumber\\
        &&-4k_1\lambda (1+f^2_0)^{\lambda-1}|f_0|\cdot|\nabla_b f|^2\cdot t^{2\lambda-1} 
        \nonumber\\
        &&-8\lambda (1+f^2_0)^{\lambda-1}|f_0|\cdot |f_\alpha A_{\bar{\beta}\bar{\alpha},\beta} |\cdot t^{2\lambda-1} 
        \nonumber\\
        &&-8\lambda (1+f^2_0)^{\lambda-1}|f_0|\cdot |f_{\alpha\beta}A_{\bar{\alpha}\bar{\beta}}|\cdot t^{2\lambda-1}  \}
        \nonumber\\
        &&-2\langle \nabla_b f,\nabla_b \mathcal{F}\rangle .
     \end{eqnarray}
     Then we estimate certain terms in (3.8) as follows
     \begin{align}
       8 \langle  \nabla&_b f_0, J\nabla_b f \rangle+2\lambda(2\lambda-1) (1+f^2_0)^{\lambda-1}|\nabla_b f_0|^2\cdot t^{2\lambda-1} 
        \nonumber\\
        &\geq-\frac{8}{\lambda(2\lambda-1)}(1+f_0^2)^{1-\lambda}|\nabla_bf |^2 \cdot t^{1-2\lambda} ,
  \\
       -8\lambda& (1+f^2_0)^{\lambda-1}|f_0|\cdot |f_\alpha A_{\bar{\beta}\bar{\alpha},\beta}|\cdot t^{2\lambda-1}  \nonumber\\
       &\geq-2k_1|\nabla_b f|^2 -2k_1\lambda^2(1+f_0^2)^{2\lambda-2}|f_0^2|\cdot t^{4\lambda-2} ,
\\
        -8\lambda& (1+f^2_0)^{\lambda-1}|f_0|\cdot |f_{\alpha\beta}A_{\bar{\alpha}\bar{\beta}} |\cdot t^{2\lambda-1} +4f_{\alpha\beta}f_{\bar{\alpha}\bar{\beta}}
        \nonumber\\
        &\geq -2k_1^2\lambda^2(1+f_0^2)^{2\lambda-2}|f_0|^2 \cdot t^{4\lambda-2}.
    \end{align}
  Finally we can get (3.6) from (3.8)-(3.11).
   \par 
   From (3.5) and the condition, we have
    \begin{eqnarray}
        (\Delta_b-\partial_t)\mathcal{G}&\geq &-\frac{\mathcal{G}}{t}+t\{ \frac{1}{m}(\Delta_b f)^2+4m f^2_0+4f_{\alpha\beta}f_{\Bar{\alpha}\Bar{\beta}}
        \nonumber\\
        &&+8 \langle \nabla_b f_0, J\nabla_b f \rangle -2k|\nabla_b f|^2
        \nonumber\\
       &&+2\lambda(2\lambda-1) (1+f^2_0)^{\lambda-1}|\nabla_b f_0|^2
        \nonumber\\
        &&-4k_1\lambda (1+f^2_0)^{\lambda-1}|f_0|\cdot|\nabla_b f|^2
        \nonumber\\
        &&-8\lambda (1+f^2_0)^{\lambda-1}|f_0|\cdot |f_\alpha A_{\bar{\beta}\bar{\alpha},\beta} |
        \nonumber\\
        &&-8\lambda (1+f^2_0)^{\lambda-1}|f_0|\cdot |f_{\alpha\beta}A_{\bar{\alpha}\bar{\beta}}| \}
        \nonumber\\
        &&-2\langle \nabla_b f,\nabla_b \mathcal{G}\rangle .
     \end{eqnarray}
     Next we are going to estimate certain terms that appear in (3.12). It is easy to prove that
     \begin{align}
        8 \langle \nabla&_b f_0, J\nabla_b f \rangle+2\lambda(2\lambda-1) (1+f^2_0)^{\lambda-1}|\nabla_b f_0|^2
        \nonumber\\
        &\geq -\frac{8}{\lambda(2\lambda-1)}(1+f_0^2)^{1-\lambda}|\nabla_bf |^2,
        \\
       -8\lambda& (1+f^2_0)^{\lambda-1}|f_0|\cdot |f_\alpha A_{\bar{\beta}\bar{\alpha},\beta}|  \nonumber\\
       &\geq-2k_1|\nabla_b f|^2-2k_1\lambda^2(1+f_0^2)^{2\lambda-2}|f_0^2|,
       \\
        -8\lambda& (1+f^2_0)^{\lambda-1}|f_0|\cdot |f_{\alpha\beta}A_{\bar{\alpha}\bar{\beta}} | +4f_{\alpha\beta}f_{\bar{\alpha}\bar{\beta}}
        \nonumber
        \\
        &\geq -2k_1^2\lambda^2(1+f_0^2)^{2\lambda-2}|f_0|^2.
    \end{align}
    Then (3.7) can be obtained by the above estimates.
    {\qed}

     ~\\
    Choose a cut-off function $\varphi \in C^{\infty}([0,\infty))$ such that
    \begin{eqnarray}
        \varphi|_{[0,1]}=1,\ \varphi|_{[2,\infty)}=0,\ -C_1^{'}|\varphi|^{\frac{1}{2}}\leq \varphi^{'}\leq 0,\ \varphi^{''}\geq -C_1^{'}. \nonumber
     \end{eqnarray}
     Set $g=\varphi^{\frac{1}{1-\mu}}$ where $\mu<1$. Direct calculations show that
     \begin{eqnarray}
        g^{'}&=&\frac{1}{1-\mu}\varphi^{\frac{\mu}{1-\mu}}\varphi^{'}=\frac{1}{1-\mu}g^\mu \varphi^{'} ,\nonumber\\
        g^{''}&=&\frac{\mu}{(1-\mu)^2}g^{2\mu-1}\varphi^{'2}+\frac{1}{1-\mu}g^\mu \varphi^{''}.
     \end{eqnarray}
    Let $r$ be the Riemannian distance and $B_p(R)$ denotes the Riemannian ball of radius $R$ centered at $p$. Put
    \begin{equation}
        \phi=g(\frac{r}{R}). \nonumber
    \end{equation}
    Assuming $R\geq 1$ and using a comparison theorem in \cite{CDRZ}, we find that
     \begin{eqnarray}
        \frac{|\nabla_b \phi|^2}{\phi^{2\mu}}&=&\frac{|\phi^{'}|^2|\nabla_b r|^2}{\phi^{2\mu} R^2}\leq \frac{C^{'}_2}{R^2} ,
        \nonumber\\
        \frac{\Delta_b \phi}{\phi^{2\mu-1}}&=&\frac{g^{''}|\nabla_b r|^2}{\phi^{2\mu-1}R^2}+\frac{g^{'}\Delta_b r}{\phi^{2\mu-1}R}\geq -\frac{C^{'}_2}{R} ,\nonumber
     \end{eqnarray}
     where $C^{'}_2$ is a constant depending on $k,k_1,\mu$. Let $\mu =3\lambda-1$, where $\frac{1}{2}<\lambda<\frac{2}{3}$, then
     \begin{eqnarray}
        \frac{|\nabla_b \phi|^2}{\phi^{6\lambda-2}}\leq \frac{C_1}{R^2},\ \ 
        \frac{\Delta_b \phi}{\phi^{6\lambda-3}}\geq -\frac{C_1}{R} ,\nonumber
     \end{eqnarray}
    where $C_1$ is a constant depending on $k,k_1,\lambda$.
    
    ~ \\
    $\mathbf{Lemma\ 3.2}$ Let $(M^{2m+1},\theta)$ be a complete noncompact pseudo-Hermitian manifold with
     \begin{equation}
                Ric_b +2(m-2)Tor_b\geq -k\ \  and\ \   |A|,|\nabla_b A|\leq k_1.\nonumber
     \end{equation}
    Let $\phi$ be defined as above with $R\geq 1$. If $\phi(x)\neq 0 $ and $\frac{1}{2}<\lambda<\frac{2}{3}$, then at $x$, we have    
    \begin{eqnarray}
         (\Delta_b-\partial_t)\phi \mathcal{F}
        &\geq& 2\langle \nabla_b (\phi \mathcal{F}),\nabla_b \phi\rangle \phi^{-1}-\frac{3C_1}{R}\phi^{6\lambda-3}\mathcal{F}-\frac{\phi \mathcal{F}_1}{t}-\frac{ 2\lambda\phi \mathcal{F}_2}{t}
         \nonumber\\
         &&-2\langle \nabla_b f,\nabla_b (\phi \mathcal{F})\rangle+2\langle \nabla_b f,\nabla_b \phi\rangle \mathcal{F}
         \nonumber\\
         &&+\phi t\{ \frac{1}{m}(\Delta_b f)^2+4mf^2_0 -(2k+2k_1)|\nabla_b f|^2
        \nonumber\\
        &&-\frac{8}{\lambda(2\lambda-1)}(1+f_0^2)^{1-\lambda}|\nabla_b f |^2\cdot t^{1-2\lambda} 
        \nonumber\\
        &&-4k_1\lambda (1+f^2_0)^{\lambda-1}|f_0|\cdot|\nabla_b f|^2 \cdot t^{2\lambda-1}
        \nonumber\\
        &&-2(k_1+k_1^2)\lambda^2(1+f_0^2)^{2\lambda-2}|f_0^2|\cdot t^{4\lambda-2} \}
        \nonumber
     \end{eqnarray}
     for $0\leq t\leq 1 $ and
     \begin{eqnarray}
         (\Delta_b-\partial_t)\phi \mathcal{G}
        &\geq& 2\langle \nabla_b (\phi \mathcal{G}),\nabla_b \phi\rangle \phi^{-1}-\frac{3C_1}{R}\phi^{6\lambda-3}\mathcal{G}
         \nonumber\\
         &&-\frac{\phi \mathcal{G}}{t}-2\langle \nabla_b f,\nabla_b (\phi \mathcal{G})\rangle+2\langle \nabla_b f,\nabla_b \phi\rangle \mathcal{G}
         \nonumber\\
         &&+\phi t\{ \frac{1}{m}(\Delta_b f)^2+4mf^2_0 -(2k+2k_1)|\nabla_b f|^2
        \nonumber\\
        &&-\frac{8}{\lambda(2\lambda-1)}(1+f_0^2)^{1-\lambda}|\nabla_b f |^2
        \nonumber\\
        &&-4k_1\lambda (1+f^2_0)^{\lambda-1}|f_0|\cdot|\nabla_b  f|^2 
        \nonumber\\
        &&-2(k_1+k_1^2)\lambda^2(1+f_0^2)^{2\lambda-2}|f_0^2| \}
        \nonumber
     \end{eqnarray}
    for $ t\geq 1 $.

    $\mathbf{Proof}$ From Lemma 3.1 and the properties of $\phi$, we have
    
    \begin{eqnarray}
         (\Delta_b-\partial_t)\phi \mathcal{F}&=&(\Delta_b \phi) \mathcal{F}+2\langle \nabla_b \phi,\nabla_b \mathcal{F} \rangle +\phi(\Delta_b-\partial_t)\mathcal{F}
         \nonumber\\
         &\geq& (\Delta_b \phi) \mathcal{F}+2\langle \nabla_b (\phi \mathcal{F}),\nabla_b \phi\rangle \phi^{-1}-\frac{2|\nabla_b \phi|^2}{\phi}\mathcal{F}
         \nonumber\\
         &&-\frac{\phi \mathcal{F}_1}{t}-\frac{2\lambda \phi \mathcal{F}_2}{t}-2\langle \nabla_b f,\nabla_b (\phi \mathcal{F})\rangle+2\langle \nabla_b f,\nabla_b \phi\rangle \mathcal{F}
         \nonumber\\
         &&+\phi t\{ \frac{1}{m}(\Delta_b f)^2+4mf^2_0 -(2k+2k_1)|\nabla_b f|^2
        \nonumber\\
        &&-\frac{8}{\lambda(2\lambda-1)}(1+f_0^2)^{1-\lambda}|\nabla_b f |^2\cdot t^{1-2\lambda} 
        \nonumber\\
        &&-4k_1\lambda (1+f^2_0)^{\lambda-1}|f_0|\cdot|\nabla_b  f|^2 \cdot t^{2\lambda-1}
        \nonumber\\
        &&-2(k_1+k_1^2)\lambda^2(1+f_0^2)^{2\lambda-2}|f_0^2|\cdot t^{4\lambda-2} \}
        \nonumber\\
        &\geq& 2\langle \nabla_b (\phi \mathcal{F}),\nabla_b \phi\rangle \phi^{-1}-\frac{3C_1}{R}\phi^{6\lambda-3}\mathcal{F}
         \nonumber\\
         &&-\frac{\phi \mathcal{F}_1}{t}-\frac{ 2\lambda \phi \mathcal{F}_2}{t}-2\langle \nabla_b f,\nabla_b (\phi \mathcal{F})\rangle+2\langle \nabla_b f,\nabla_b \phi\rangle \mathcal{F}
         \nonumber\\
         &&+\phi t\{ \frac{1}{m}(\Delta_b f)^2+4mf^2_0 -(2k+2k_1)|\nabla_b f|^2
        \nonumber\\
        &&-\frac{8}{\lambda(2\lambda-1)}(1+f_0^2)^{1-\lambda}|\nabla_b f |^2\cdot t^{1-2\lambda} 
        \nonumber\\
        &&-4k_1\lambda (1+f^2_0)^{\lambda-1}|f_0|\cdot|\nabla_b  f|^2 \cdot t^{2\lambda-1}
        \nonumber\\
        &&-2(k_1+k_1^2)\lambda^2(1+f_0^2)^{2\lambda-2}|f_0^2|\cdot t^{4\lambda-2} \
        \nonumber
     \end{eqnarray}
     for $0\leq t\leq 1$ and 
    \begin{eqnarray}
         (\Delta_b-\partial_t)\phi \mathcal{G}&=&(\Delta_b \phi) \mathcal{G}+2\langle \nabla_b \phi,\nabla_b \mathcal{G} \rangle +\phi(\Delta_b-\partial_t)\mathcal{G} 
         \nonumber\\
         &\geq& (\Delta_b \phi) \mathcal{G}+2\langle \nabla_b (\phi \mathcal{G}),\nabla_b \phi\rangle \phi^{-1}-\frac{2|\nabla_b \phi|^2}{\phi}\mathcal{G}
         \nonumber\\
         &&-\frac{\phi \mathcal{G}}{t}-2\langle \nabla_b f,\nabla_b (\phi \mathcal{G})\rangle+2\langle \nabla_b f,\nabla_b \phi\rangle \mathcal{G}
         \nonumber\\
         &&+\phi t\{ \frac{1}{m}(\Delta_b f)^2+4mf^2_0 -(2k+2k_1)|\nabla_b f|^2
        \nonumber\\
        &&-\frac{8}{\lambda(2\lambda-1)}(1+f_0^2)^{1-\lambda}|\nabla_b f |^2 
        \nonumber\\
        &&-4k_1\lambda (1+f^2_0)^{\lambda-1}|f_0|\cdot|\nabla_b  f|^2 
        \nonumber\\
        &&-2(k_1+k_1^2)\lambda^2(1+f_0^2)^{2\lambda-2}|f_0^2|\}
        \nonumber\\
        &\geq& 2\langle \nabla_b (\phi \mathcal{G}),\nabla_b \phi\rangle \phi^{-1}-\frac{3C_1}{R}\phi^{6\lambda-3}\mathcal{G}
         \nonumber\\
         &&-\frac{\phi \mathcal{G}}{t}-2\langle \nabla_b f,\nabla_b (\phi \mathcal{G})\rangle+2\langle \nabla_b f,\nabla_b \phi\rangle \mathcal{G}
         \nonumber\\
         &&+\phi t\{ \frac{1}{m}(\Delta_b f)^2+4mf^2_0 -(2k+2k_1)|\nabla_b f|^2
        \nonumber\\
        &&-\frac{8}{\lambda(2\lambda-1)}(1+f_0^2)^{1-\lambda}|\nabla_b f |^2
        \nonumber\\
        &&-4k_1\lambda (1+f^2_0)^{\lambda-1}|f_0|\cdot|\nabla_b  f|^2 
        \nonumber\\
        &&-2(k_1+k_1^2)\lambda^2(1+f_0^2)^{2\lambda-2}|f_0^2| \}
        \nonumber
     \end{eqnarray}
    for $t\geq 1 $. This completes the proof.
     {\qed }
     ~\par
     Now we are ready to consider the first case of the gradient estimate, that is, $0<t\leq 1$.
     
     ~\\
    $\mathbf{Proposition\ 3.3}$ Let $(M^{2m+1},HM,J,\theta)$ be a complete noncompact pseudo-Hermitian manifold with
    \begin{eqnarray}
        Ric_b +2(m-2)Tor_b\geq -k\ \ and\ \  |A|,|\nabla_b A|\leq k_1.\nonumber
     \end{eqnarray}
     and $u$ be a positive solution of the CR heat equation
    \begin{eqnarray}
        \frac{\partial u}{\partial t}=\Delta_b u  \nonumber
     \end{eqnarray}
    on $B_p(2R)\times (0,1]$ with $R\geq 1$. Then for any constant $\frac{1}{2}<\lambda<\frac{2}{3}$ and any constant $\delta>1+\frac{4}{m\lambda(2\lambda-1)}$, there exists a constant $C_3^{'}$ depends on $m,k,k_1,\lambda,\delta$, such that 
    \begin{eqnarray}
      \frac{|\nabla_b u|^2}{u^2}+t^{2\lambda-1}(1+\frac{u^2_0}{u^2})^\lambda-\delta \frac{u_t}{u}\leq \frac{C_3^{'}}{t}(1+\frac{1}{R^\lambda})
    \end{eqnarray}
     on $B_p(R)\times (0,1]$.
     
     ~\\
$\mathbf{Proof}$ 
     Let $(x_1,t_1)$ be the maximum point of $\phi \mathcal{F}$ on $B_p(2R)\times [0,1]$. Without loss of generality, we may assume that $(\phi \mathcal{F})(x_1,t_1)>0$, otherwise the conclusion follows trivially. At $(x_1,t_1)$, we have $\nabla (\phi \mathcal{F})=0,\ \partial_t (\phi\mathcal{F})\geq 0$ and $\Delta_b (\phi \mathcal{F})\leq 0.$ Using Lemma
    3.2 and evaluating the inequality at $(x_1,t_1)$, we obtain
     \begin{eqnarray}
       0&\geq& -\frac{3C_1}{R}\phi^{6\lambda-3}\mathcal{F}
         -\frac{\phi \mathcal{F}_1}{t_1}-\frac{2\lambda \phi \mathcal{F}_2}{t_1}+2\langle \nabla_b f,\nabla_b \phi\rangle \mathcal{F}
         \nonumber\\
         &&+\phi t_1\{\frac{1}{m}(\Delta_b f)^2+4mf^2_0 
        \nonumber\\
        &&-\frac{8}{\lambda(2\lambda-1)}(1+f_0^2)^{1-\lambda}|\nabla_b f |^2\cdot t_1^{1-2\lambda} -2k|\nabla_b  f|^2
        \nonumber\\
        &&-2k_1|\nabla_b  f|^2-4k_1\lambda (1+f^2_0)^{\lambda-1}|f_0|\cdot|\nabla_b  f|^2 \cdot t_1^{2\lambda-1}
        \nonumber\\
        &&-2(k_1+k_1^2)\lambda^2(1+f_0^2)^{2\lambda-2}|f_0^2|\cdot t_1^{4\lambda-2} \}.
        \nonumber
     \end{eqnarray}
     Multiplying $\phi t_1$ yields  
     \begin{eqnarray}
       0&\geq& -\frac{3C_1}{R}\phi^{6\lambda-2}t_1 \mathcal{F}
         -\phi^2 \mathcal{F}_1-2\lambda \phi^2 \mathcal{F}_2-\frac{2C_1}{R}|\nabla_b f|\cdot \phi^{3\lambda}t_1 \mathcal{F} 
         \nonumber\\
         &&+\phi^2 t^2_1\{ \frac{1}{m}(\Delta_b f)^2+4mf^2_0 
        \nonumber\\
        &&-\frac{8}{\lambda(2\lambda-1)}(1+f_0^2)^{1-\lambda}|\nabla_b f |^2\cdot t_1^{1-2\lambda} -2k|\nabla_b  f|^2
        \nonumber\\
        &&-2k_1|\nabla_b f|^2-4k_1\lambda (1+f^2_0)^{\lambda-1}|f_0|\cdot|\nabla_b  f|^2 \cdot t_1^{2\lambda-1}
        \nonumber\\
        &&-2(k_1+k_1^2)\lambda^2(1+f_0^2)^{2\lambda-2}|f_0^2|\cdot t_1^{4\lambda-2} \}
       \\
        &\geq& -\frac{3C_1}{R}\phi^{6\lambda-2}t_1 \mathcal{F}
         -\phi^2 \mathcal{F}_1-2\lambda \phi^2 \mathcal{F}_2-\frac{2C_1}{R}|\nabla_b f|\cdot \phi^{3\lambda} t_1 \mathcal{F} 
         \nonumber\\
         &&+\phi^2 t^2_1\{ \frac{1}{m}(\Delta_b f)^2+4mf^2_0 -2\epsilon |\nabla_b  f|^4
        \nonumber\\
        &&-\frac{16}{\epsilon\lambda^2(2\lambda-1)^2}(1+f_0^2)^{2-2\lambda}t_1^{2-4\lambda} -(2k+2k_1 )|\nabla_b  f|^2
        \nonumber\\
        &&-4\epsilon^{-1}k^2_1\lambda^2 (1+f^2_0)^{2\lambda-2}|f_0|^2\cdot t_1^{4\lambda-2}
        \nonumber\\
        &&-2(k_1+k_1^2)\lambda^2(1+f_0^2)^{2\lambda-2}|f_0^2|\cdot t_1^{4\lambda-2}\},
     \end{eqnarray}
     since 
     \begin{eqnarray}
        && -\frac{8}{\lambda(2\lambda-1)}(1+f_0^2)^{1-\lambda}|\nabla_b f |^2  \cdot t_1^{1-2\lambda}
        \nonumber\\
        &\geq& -\epsilon |\nabla_b  f|^4-\frac{16}{\epsilon\lambda^2(2\lambda-1)^2}(1+f_0^2)^{2-2\lambda}\cdot t_1^{2-4\lambda}
     \end{eqnarray}
     \begin{eqnarray}
       &&-4k_1\lambda (1+f^2_0)^{\lambda-1}|f_0|\cdot|\nabla_b  f|^2 \cdot t_1^{2\lambda-1}
       \nonumber\\
       &\geq &-\epsilon |\nabla_b f|^4-4\epsilon^{-1} \lambda^2 k_1^2(1+f_0^2)^{2\lambda-2}f_0^2\cdot t_1^{4\lambda-2},
    \end{eqnarray}
    where $\epsilon$ is a constant to be determined. Hence we have 
     \begin{eqnarray}
       0
        &\geq& -\frac{3C_1}{R}\phi^{6\lambda-2}t_1 \mathcal{F}
         -\phi^2 \mathcal{F}_1-2\lambda \phi^2 \mathcal{F}_2-\frac{2C_1}{R}|\nabla_b f|\cdot \phi^{3\lambda} t_1 \mathcal{F} 
         \nonumber\\
         &&+\phi^2 t^2_1\{ \frac{1}{m}(|\nabla_b f|^2-f_t)^2+4mf^2_0 -2\epsilon |\nabla_b  f|^4
        \nonumber\\
        &&-\frac{16}{\epsilon\lambda^2(2\lambda-1)^2}(1+f_0^2)^{2-2\lambda}t_1^{2-4\lambda} -(2k+2k_1 )|\nabla_b  f|^2
        \nonumber\\
        &&-4\epsilon^{-1}k^2_1\lambda^2 (1+f^2_0)^{2\lambda-2}|f_0|^2\cdot t_1^{4\lambda-2}
        \nonumber\\
        &&-2(k_1+k_1^2)\lambda^2(1+f_0^2)^{2\lambda-2}|f_0^2|\cdot t_1^{4\lambda-2}\}.
     \end{eqnarray}
     Note that it is difficult to estimate $\mathcal{F}$ directly. Let us recall the method
    of Cao-Yau \cite{CaY} for weakly elliptic operators on closed manifolds. Translating their idea to pseudo-Hermitian case, they actually tried to control $\mathcal{F}$ by either $\delta_0|\nabla_b f|^2-\delta f_t$ ($\delta>\delta_0>1$) or $t_1^{2\lambda-1}(1+f^2_0)^\lambda$. However, Cao-Yau’s estimates
cannot be applied directly to give the required inequalities in the complete
noncompact case. We have to treat some extra terms appearing in (3.22), e.g., $\frac{2C_1}{R}|\nabla_b f|\cdot \phi^{3\lambda} t_1 \mathcal{F} $. Our discussion will be divided into two cases according to the sign of $f_t$ as follows.
     
    ~\\
     $\mathbf{Case\  A}$ $f_t< 0$ at the maximum point $(x_1,t_1)$.
     
     ~\\
    In this case, we have 
    \begin{eqnarray}
       (|\nabla_b f|^2-f_t)^2\geq |\nabla_b f|^4+f^2_t. \nonumber
    \end{eqnarray}
     Then (3.22) becomes 
     \begin{eqnarray}
       0&\geq&-\frac{3C_1}{R}\phi^{6\lambda-2}t_1 \mathcal{F}
       -\phi^2 \mathcal{F}_1-2\lambda \phi^2 \mathcal{F}_2+2m\phi^2t_1^2f_0^2+\frac{\phi^2 t_1^2}{2m}(|\nabla_b f|^2-f_t)^2
         \nonumber\\
         &&+\phi^2 t^2_1\{ (\frac{1}{2m}-2\epsilon)|\nabla_b f|^4
        -2k|\nabla_b  f|^2-2k_1|\nabla_b  f|^2\}
        \nonumber\\
        &&+\phi^2 t^2_1\{2mf_0^2-\frac{16}{\epsilon\lambda^2(2\lambda-1)^2}(1+f_0^2)^{2-2\lambda}t_1^{2-4\lambda} 
        \nonumber\\
        &&-4\epsilon^{-1}k^2_1\lambda^2 (1+f^2_0)^{2\lambda-2}|f_0|^2\cdot t_1^{4\lambda-2}
        \nonumber\\
        &&-2(k_1+k_1^2)\lambda^2(1+f_0^2)^{2\lambda-2}|f_0^2|\cdot t_1^{4\lambda-2}\}
        \nonumber\\
        &&+\phi^2 t^2_1 \frac{1}{2m}f_t^2-\frac{2C_1}{R}|\nabla_b f|\cdot \phi^{3\lambda}t_1 \mathcal{F}.
     \end{eqnarray}
     Following the idea in \cite{CaY}, we want to control $\mathcal{F}$ by either $\delta_0|\nabla_b f|^2-\delta f_t$ ($\delta>\delta_0>1$) or $t_1^{2\lambda-1}(1+f^2_0)^\lambda$.
     
     ~\\
     $\mathbf{(A1)}$ Suppose $\delta_0|\nabla_b f|^2-\delta f_t\geq t_1^{2\lambda-1}(1+f^2_0)^\lambda$. Clearly 
      \begin{eqnarray}
       \mathcal{F}\leq t_1(|\nabla_b f|^2-\delta f_t+\delta_0|\nabla_b f|^2-\delta f_t)=t_1((\delta_0+1)|\nabla_b f|^2-2\delta f_t) .\nonumber
     \end{eqnarray}
    Consequently the last term of (3.23) can be estimate by
     \begin{eqnarray}
       &&\frac{2C_1}{R}|\nabla_b f|\cdot \phi^{3\lambda}t_1 \mathcal{F} 
       \nonumber\\
       &\leq &\frac{2C_1}{R}(\delta_0+1)t_1^2 (\phi |\nabla_b f|^2)^{\frac{3}{2}}+\frac{4C_1}{R}t_1^2 \delta(\phi|\nabla_b f|^2)^{\frac{1}{2}}\cdot |\phi f_t| \nonumber\\
       &\leq& \frac{2C_1}{R}(\delta_0+1)t_1^2 (\phi |\nabla_b f|^2)^{\frac{3}{2}}+\phi^2 t^2_1 \frac{1}{2m}f_t^2+\frac{8mC_1^2}{R^2}\delta^2 t_1^2 \phi |\nabla_b f|^2.
    \end{eqnarray}
  Noting that $0<t_1\leq 1 $, (3.23) and (3.24) yield that 
    \begin{eqnarray}
       0&\geq&-\frac{3C_1}{R}\phi^{6\lambda-2}t_1 \mathcal{F}
         -\phi^2 \mathcal{F}_1-2\lambda \phi^2 \mathcal{F}_2+2m\phi^2t_1^2f_0^2+\frac{\phi^2 t_1^2}{2m}(|\nabla_b f|^2-f_t)^2
         \nonumber\\
         &&+ t^2_1\{ (\frac{1}{2m}-2\epsilon)(\phi|\nabla_b  f|^2)^2
        -2(k+k_1+\frac{4mC_1^2}{R^2}\delta^2)\phi|\nabla_b  f|^2
        \nonumber\\
        &&-\frac{2C_1}{R}(\delta_0+1) (\phi |\nabla_b f|^2)^{\frac{3}{2}}\}
        \nonumber\\
        &&+\phi^2 \{2m(t_1f_0)^2-\frac{16}{\epsilon\lambda^2(2\lambda-1)^2}(1+f_0^2)^{2-2\lambda}t_1^{4-4\lambda} 
        \nonumber\\
        &&-4\epsilon^{-1}k^2_1\lambda^2 (1+f^2_0)^{2\lambda-2}|f_0|^2\cdot t_1^{4\lambda}
        \nonumber\\
        &&-2(k_1+k_1^2)\lambda^2(1+f_0^2)^{2\lambda-2}|f_0^2|\cdot t_1^{4\lambda}\}.
    \end{eqnarray}
    Let $\epsilon<\frac{1}{4m}$. By (3.25), we have 
    \begin{eqnarray}
       0&\geq&-\frac{3C_1}{R}\phi^{6\lambda-2}t_1 \mathcal{F}
         -\phi^2 \mathcal{F}_1-2\lambda \phi ^2 \mathcal{F}_2+\frac{\phi^2 t_1^2}{2m}(|\nabla_b f|^2-f_t)^2-C_2
         \nonumber\\
         &&-t^2_1\{ C_2+\frac{C_2}{R^2}+\frac{C_2}{R^4}\},
    \end{eqnarray}
    where $C_2$ is a constant depending on $m,k,k_1,\delta,\delta_0,\lambda$. Let $x=\phi (\delta_0 |\nabla_b  f|^2-\delta f_t)(x_1,t_1)$, then 
    \begin{eqnarray}
      \phi \mathcal{F}_1\leq t_1 x,\ \phi \mathcal{F}_2\leq t_1 x,\  \phi(|\nabla_b f|^2-f_t)\geq \frac{1}{\delta} x.\nonumber
    \end{eqnarray}
    Therefore we find that 
    \begin{eqnarray}
      0\geq \frac{t_1^2}{2m\delta^2}x^2-(2\lambda+1)t_1x-\frac{6C_1}{R}t_1^2 x-t_1^2(C_2+\frac{C_2}{R^2}+\frac{C_2}{R^4})-C_2 ,
    \end{eqnarray}
    which implies that 
    \begin{eqnarray}
      t_1 x \leq C_3(1+\frac{1}{R}) ,\nonumber
    \end{eqnarray}
    and 
    \begin{eqnarray}
      \phi \mathcal{F}\leq 2t_1 x \leq 2C_3(1+\frac{1}{R}) ,
    \end{eqnarray}
     where $C_3$ is a constant depending on $m,k,k_1,\lambda,\delta_0,\delta$.
     
     ~\\
  $\mathbf{(A2)}$ Suppose $\delta_0|\nabla_b f|^2-\delta f_t\leq t_1^{2\lambda-1}(1+f^2_0)^\lambda$. Then
    \begin{eqnarray}
      \phi \mathcal{F}&=& \phi t_1(\delta_0|\nabla_b f|^2+t_1^{2\lambda-1}(1+f_0^2)^\lambda-\delta f_t+(1-\delta_0)|\nabla_b f|^2)
      \nonumber\\
      &\leq& 2\phi t_1^{2\lambda}(1+f_0^2)^\lambda. \nonumber
    \end{eqnarray}
    Since $f_t< 0,t_1\leq1 $ and $\lambda>\frac{1}{2}$, the assumption (A2) implies
    \begin{eqnarray}
       (1+f_0^2)^\lambda\geq \delta_0 |\nabla_b  f|^2.\nonumber
    \end{eqnarray}
   The last term of (3.23) is bounded by
     \begin{eqnarray}
      \frac{2C_1}{R}|\nabla_b f|\cdot \phi^{3\lambda}t_1 \mathcal{F} 
       &\leq &\frac{4C_1}{R}\phi^{3\lambda}t^2_1|\nabla_b  f|\cdot (1+f_0^2)^\lambda
       \nonumber\\
       &\leq&\frac{4C_1}{R\sqrt{\delta_0}}\phi^{3\lambda}t^2_1\cdot (1+f_0^2)^{\frac{3}{2}\lambda}
       \nonumber\\
       &=&\frac{4C_1}{R\sqrt{\delta_0}}t^2_1\cdot [\phi^2(1+f_0^2)]^{\frac{3}{2}\lambda}.
       \nonumber
    \end{eqnarray}
    Hence (3.23) becomes 
     \begin{eqnarray}
       0&\geq&-\frac{3C_1}{R}\phi^{6\lambda-2}t_1 \mathcal{F}
         -\phi^2 \mathcal{F}_1-2\lambda \phi^2 \mathcal{F}_2+2m\phi^2t_1^2f_0^2
         \nonumber\\
         &&+ t^2_1\{ (\frac{1}{2m}-2\epsilon)(\phi|\nabla_b  f|^2)^2
        -2(k+k_1)\phi|\nabla_b  f|^2 \}
        \nonumber\\
        &&+ \{2m(\phi t_1f_0)^2-\frac{16}{\epsilon\lambda^2(2\lambda-1)^2}(1+f_0^2)^{2-2\lambda}t_1^{4-4\lambda} \phi^2
        \nonumber\\
        &&-4\epsilon^{-1}k^2_1\lambda^2 (1+f^2_0)^{2\lambda-2}|f_0|^2\cdot t_1^{4\lambda}\phi^2
        \nonumber\\
        &&-2(k_1+k_1^2)\lambda^2(1+f_0^2)^{2\lambda-2}|f_0^2|\cdot t_1^{4\lambda}\phi^2
        \nonumber\\
        &&-\frac{4C_1}{R\sqrt{\delta_0}}t^2_1\cdot [\phi^2(1+f_0^2)]^{\frac{3}{2}\lambda}\}.\
    \end{eqnarray}
    Letting $\epsilon<\frac{1}{4m}$ and noting that  $\frac{1}{2}<\lambda <\frac{2}{3}$, we have
    \begin{eqnarray}
       0&\geq&-\frac{3C_1}{R}\phi^{6\lambda-2}t_1 \mathcal{F}
         -\phi^2 \mathcal{F}_1-2\lambda\phi^2 \mathcal{F}_2 +2m\phi^2t_1^2f_0^2-C_4(1+\frac{1}{R}),\nonumber
    \end{eqnarray}
    where $C_4$ is a constant depending on $m,k,k_1,\delta_0,\lambda$. Set $y=\phi |f_0|$, then
    \begin{eqnarray}
       0\geq2mt_1^2 y^2-(2\lambda+1+\frac{6C_1}{R}) (t_1 y)^{2\lambda}- C_5(1+\frac{1}{R}+\frac{1}{R^4}),
    \end{eqnarray}
    which yields that 
    \begin{eqnarray}
       t_1 y\leq& C_6(1+\frac{1}{\sqrt{R}}) ,\nonumber
    \end{eqnarray}
    hence 
    \begin{eqnarray}
       \phi \mathcal{F}\leq 2\phi t_1^{2\lambda}(1+f^2_0)^\lambda\leq C_7(1+\frac{1}{R^\lambda}) ,
    \end{eqnarray}
    where $C_5,C_6,C_7$ are constants depending on $m,k,k_1,\delta_0,\lambda$.
    
    ~\\
    $\mathbf{Case\ B}$ $f_t\geq 0$ at the maximal point $(x_1,t_1)$.
    
    ~\\
    $\mathbf{(B1)} $ Suppose $\delta_0|\nabla_b f|^2-\delta f_t\geq 0$. We observe that 
    \begin{eqnarray}
       (|\nabla_b  f|^2-f_t)^2&=&\{ \frac{1}{\delta}(\delta_0|\nabla_b  f|^2-\delta f_t)+(1-\frac{\delta_0}{\delta})|\nabla_b  f|^2 \}^2 \nonumber\\
       &\geq& \frac{1}{\delta^2}(\delta_0|\nabla_b f|^2-\delta f_t)^2+(1-\frac{\delta_0}{\delta})^2|\nabla_b f|^4.
    \end{eqnarray}
  Thus (3.22) becomes 
    \begin{eqnarray}
       0 &\geq& -\frac{3C_1}{R}\phi^{6\lambda-2}t_1 \mathcal{F}
         -\phi^2 \mathcal{F}_1-2\lambda \phi^2 \mathcal{F}_2+\frac{\phi^2t_1^2}{m\delta^2}(\delta_0|\nabla_b  f|^2-\delta f_t)^2
         \nonumber\\
         &&+\phi^2 t^2_1\{ 
        (\frac{(\delta-\delta_0)^2}{m\delta^2}-2\epsilon )|\nabla_b  f|^4-2(k+k_1)|\nabla_b f|^2\}
        \nonumber\\
        &&+\phi^2t_1^2\{4mf^2_0 -\frac{16}{\epsilon\lambda^2(2\lambda-1)^2}(1+f_0^2)^{2-2\lambda}t_1^{2-4\lambda} 
        \nonumber\\
        &&-4\epsilon^{-1}k^2_1\lambda^2 (1+f^2_0)^{2\lambda-2}|f_0|^2\cdot t_1^{4\lambda-2}
        \nonumber\\
        &&-2(k_1+k_1^2)\lambda^2(1+f_0^2)^{2\lambda-2}|f_0^2|\cdot t_1^{4\lambda-2}\}
        \nonumber\\
        &&-\frac{2C_1}{R}|\nabla_b f|\cdot \phi^{3\lambda} t_1 \mathcal{F} .
     \end{eqnarray}
       In this case, the proof is almost the same as that for the case of  $f_t<0$.  
       
     ~\\
     $\mathbf{(B1}$-$\mathbf{1)} $ Suppose $\delta_0|\nabla_b f|^2-\delta f_t\geq t_1^{2\lambda-1}(1+f^2_0)^\lambda$. Using the assumption that $f_t\geq 0$, we have
      \begin{eqnarray}
       \mathcal{F}\leq t_1(|\nabla_b f|^2-\delta f_t+\delta_0|\nabla_b f|^2-\delta f_t)=t_1((\delta_0+1)|\nabla_b f|^2) .\nonumber
     \end{eqnarray}
     Hence the last term of (3.33) can be estimated by
     \begin{eqnarray}
       \frac{2C_1}{R}|\nabla_b f|\cdot \phi^{3\lambda}t_1 \mathcal{F} 
       \leq \frac{2C_1}{R}(\delta_0+1)t_1^2 (\phi |\nabla_b f|^2)^{\frac{3}{2}}. 
    \end{eqnarray}
   Therefore (3.33) and (3.34) yield that 
     \begin{eqnarray}
       0 &\geq& -\frac{3C_1}{R}\phi^{6\lambda-2}t_1 \mathcal{F}
         -\phi^2 \mathcal{F}_1-2\lambda \phi^2 \mathcal{F}_2+\frac{\phi^2t_1^2}{m\delta^2}(\delta_0|\nabla_b  f|^2-\delta f_t)^2
         \nonumber\\
         &&+ t^2_1\{ 
        (\frac{(\delta-\delta_0)^2}{m\delta^2}-2\epsilon )\phi^2|\nabla_b  f|^4-2(k+k_1)\phi^2|\nabla_b f|^2
        \nonumber\\
        &&-\frac{2C_1}{R}(\delta_0+1) (\phi |\nabla_b f|^2)^{\frac{3}{2}}\}
        \nonumber\\
        &&+\phi^2\{4m(t_1|f_0|)^2 -\frac{16}{\epsilon\lambda^2(2\lambda-1)^2}(1+f_0^2)^{2-2\lambda}t_1^{4-4\lambda} 
        \nonumber\\
        &&-4\epsilon^{-1}k^2_1\lambda^2 (1+f^2_0)^{2\lambda-2}|f_0|^2\cdot t_1^{4\lambda}
        \nonumber\\
        &&-2(k_1+k_1^2)\lambda^2(1+f_0^2)^{2\lambda-2}|f_0^2|\cdot t_1^{4\lambda}\}.   
     \end{eqnarray}
    Let $\epsilon<\frac{(\delta-\delta_0)^2}{2m\delta^2}$. By (3.35), we have 
    \begin{eqnarray}
       0&\geq&-\frac{3C_1}{R}\phi^{6\lambda-2}t_1 \mathcal{F}
         -\phi^2 \mathcal{F}_1-2\lambda \phi ^2 \mathcal{F}_2+\frac{\phi^2 t_1^2}{m\delta^2}(\delta_0|\nabla_b f|^2-\delta f_t)^2
         \nonumber\\
         &&-t^2_1\{ C_8+\frac{C_8}{R^4}\}-C_8,
    \end{eqnarray}
    where $C_8$ is a constant depending on $m,k,k_1,\delta,\delta_0,\lambda$. Let $x=\phi (\delta_0 |\nabla_b  f|^2-\delta f_t)(x_1,t_1)$, then 
    \begin{eqnarray}
      \phi \mathcal{F}_1\leq t_1 x,\ \phi \mathcal{F}_2\leq t_1 x,\  \phi(|\nabla_b f|^2-f_t)\geq \frac{1}{\delta} x.\nonumber
    \end{eqnarray}
    Therefore we find that 
    \begin{eqnarray}
      0\geq \frac{t_1^2}{m\delta^2}x^2-(2\lambda+1)t_1x-\frac{6C_1}{R}t_1^2 x-t_1^2(C_8+\frac{C_8}{R^4})-C_8 ,
    \end{eqnarray}
    which implies that 
    \begin{eqnarray}
      t_1 x \leq C_9(1+\frac{1}{R}) ,\nonumber
    \end{eqnarray}
    and 
    \begin{eqnarray}
      \phi \mathcal{F}\leq 2t_1 x \leq 2C_9(1+\frac{1}{R}) ,
    \end{eqnarray}
     where $C_9$ is a constant depending on $m,k,k_1,\lambda,\delta_0,\delta$.
     
     ~\\
   $\mathbf{(B1}$-$\mathbf{2)}$  Suppose $\delta_0|\nabla_b f|^2-\delta f_t\leq t_1^{2\lambda-1}(1+f^2_0)^\lambda$. Then
    \begin{eqnarray}
      \phi \mathcal{F}&=& \phi t_1(\delta_0|\nabla_b f|^2+t_1^{2\lambda-1}(1+f_0^2)^\lambda-\delta f_t+(1-\delta_0)|\nabla_b f|^2)
      \nonumber\\
      &\leq& 2\phi t_1^{2\lambda}(1+f_0^2)^\lambda. \nonumber
    \end{eqnarray}
   The last term of (3.33) can be estimated by
     \begin{eqnarray}
      \frac{2C_1}{R}|\nabla_b f|\cdot \phi^{3\lambda}t_1 \mathcal{F} 
       &\leq &\frac{4C_1}{R}\phi^{3\lambda}t^2_1|\nabla_b  f|\cdot (1+f_0^2)^\lambda
       \nonumber\\
       &\leq&\frac{C_{10}}{R}\phi^{3\lambda}t^2_1\cdot [|\nabla_b f|^3+(1+f_0^2)^{\frac{3}{2}\lambda}]
       \nonumber\\
       &\leq&\frac{C_{10}}{R}t^2_1\cdot [(\phi|\nabla_b f|^2)^{\frac{3}{2}}+|\phi^2(1+f_0^2)|^{\frac{3}{2}\lambda}],
       \nonumber
    \end{eqnarray}
    where we use the Young's inequality $ab\leq \frac{1}{3}a^3+\frac{2}{3}b^{\frac{3}{2}}(a,b\geq 0) $ in the second inequality. Hence (3.33) becomes 
     \begin{eqnarray}
       0&\geq&-\frac{3C_1}{R}\phi^{6\lambda-2}t_1 \mathcal{F}
         -\phi^2 \mathcal{F}_1-2\lambda \phi^2 \mathcal{F}_2+2m\phi^2t_1^2f_0^2
         \nonumber\\
         &&+ t^2_1\{ (\frac{(\delta-\delta_0)^2}{m\delta^2}-2\epsilon)(\phi|\nabla_b  f|^2)^2
        -2(k+k_1)\phi|\nabla_b  f|^2 
        \nonumber\\
        &&-\frac{C_{10}}{R}(\phi|\nabla_b f|^2)^{\frac{3}{2}}\}
        \nonumber\\
        &&+ \{2m(\phi t_1|f_0|)^2-\frac{16}{\epsilon\lambda^2(2\lambda-1)^2}(1+f_0^2)^{2-2\lambda}t_1^{4-4\lambda} \phi^2
        \nonumber\\
        &&-4\epsilon^{-1}k^2_1\lambda^2 (1+f^2_0)^{2\lambda-2}|f_0|^2\cdot t_1^{4\lambda}\phi^2
        \nonumber\\
        &&-2(k_1+k_1^2)\lambda^2(1+f_0^2)^{2\lambda-2}|f_0^2|\cdot t_1^{4\lambda}\phi^2
        \nonumber\\
        &&-\frac{C_{10}}{R}t^2_1\cdot [\phi^2(1+f_0^2)]^{\frac{3}{2}\lambda}\}.\
    \end{eqnarray}
   Choosing $\epsilon<\frac{(\delta-\delta_0)^2}{2m\delta^2}$ and noting that $\lambda<\frac{2}{3}$ and $t_1\leq 1$, we have 
    \begin{eqnarray}
       0\geq-\frac{3C_1}{R}\phi^{6\lambda-2}t_1 \mathcal{F}
         -\phi^2 \mathcal{F}_1-2\lambda\phi^2 \mathcal{F}_2 +2m\phi^2t_1^2f_0^2-C_{11}(1+\frac{1}{R}+\frac{1}{R^4}),
    \end{eqnarray}
    where $C_{11}$ is a constant depending on $m,k,k_1,\delta,\delta_0,\lambda$. Set $y=\phi |f_0|$. Then
    \begin{eqnarray}
       0\geq2mt_1^2 y^2-(2\lambda+1+\frac{6C_1}{R}) (t_1 y)^{2\lambda}- C_{12}(1+\frac{1}{R}+\frac{1}{R^4}), \nonumber
    \end{eqnarray}
    which yields that 
    \begin{eqnarray}
       t_1 y\leq C_{13}(1+\frac{1}{\sqrt{R}}) ,\nonumber
    \end{eqnarray}
    hence 
    \begin{eqnarray}
       \phi \mathcal{F}\leq 2\phi t_1^{2\lambda}(1+f^2_0)^\lambda\leq C_{14}(1+\frac{1}{R^\lambda}) ,
    \end{eqnarray}
    where $C_{12},C_{13},C_{14}$ are constants depending on $m,k,k_1,\delta,\delta_0,\lambda$.
    
    ~\\
   $\mathbf{(B2)}$ Suppose $\delta_0|\nabla_b f|^2-\delta f_t\leq 0$. In this case, we can assume that 
    \begin{eqnarray}
       (\delta_0-1)|\nabla_b f|^2\leq t_1^{2\lambda-1}(1+f_0^2)^\lambda.
    \end{eqnarray}
    Otherwise 
    \begin{eqnarray}
        \mathcal{F}&=&t_1(|\nabla_b f|^2+t_1^{2\lambda-1}(1+f_0^2)^\lambda-\delta f_t) \nonumber\\
        &\leq&t_1(\delta_0|\nabla_b f|^2-\delta f_t)\leq 0, \nonumber
    \end{eqnarray}
    and thus the conclusion of Proposition 3.3 follows trivially. From (3.18) and (3.42), we have 
   
    \begin{eqnarray}
        0
        &\geq& -\frac{3C_1}{R}\phi^{6\lambda-2}t_1 \mathcal{F}
         -\phi^2 \mathcal{F}_1-2\lambda\phi^2\mathcal{F}_2+2m\phi^2 t^2_1f^2_0
         \nonumber\\
         &&+\{(2m-\frac{8}{\lambda(2\lambda-1)(\delta_0-1)})\phi^2t_1^2 f^2_0 
        -\frac{8}{\lambda(2\lambda-1)(\delta_0-1)}
        \nonumber\\
        &&-(2k+2k_1)\frac{t_1^{2\lambda+1}}{\delta_0-1}\phi^2(1+f_0^2)^\lambda
        \nonumber\\
        &&-\frac{4k_1\lambda t_1^{4\lambda}}{\delta_0-1}\phi^2(1+f_0^2)^{2\lambda-1}|f_0|
        \nonumber\\
        &&-(2k_1+2k_1^2)\lambda^2\phi^2(1+f^2_0)^{2\lambda-2}f^2_0 \cdot t_1^{4\lambda}
        \nonumber\\
        &&-\frac{2C_1t_1^{3\lambda+\frac{1}{2}}}{ \sqrt{\delta_0-1} R}[\phi^2(1+f_0^2)]^{\frac{3}{2}\lambda}\}.
    \end{eqnarray}
   Choosing $\delta_0>1+\frac{4}{m\lambda(2\lambda-1)}$ and noting that $\lambda<\frac{2}{3}$, we obtain
    \begin{eqnarray}
       0\geq-\frac{3C_1}{R}\phi^{6\lambda-2}t_1 \mathcal{F}
         -\phi^2 \mathcal{F}_1-2\lambda\phi^2 \mathcal{F}_2+2m\phi^2t_1^2f_0^2-C_{15}(1+\frac{1}{R}),
    \end{eqnarray}
    where $C_{15}$ is a constant depending on $m,k,k_1,\delta_0,\lambda$. Set $y=\phi |f_0|$, we get
    \begin{eqnarray}
       0\geq2mt_1^2 y^2-(2\lambda+\frac{3C_1}{R}) (t_1y)^{2\lambda}- C_{16}(1+\frac{1}{R}),
         \nonumber
    \end{eqnarray}
    which yields that 
    \begin{eqnarray}
      t_1 y \leq C_{17}(1+\frac{1}{\sqrt{R}}) ,\nonumber
    \end{eqnarray}
    and 
    \begin{eqnarray}
       \phi \mathcal{F}\leq \phi t_1^{2\lambda}(1+f^2_0)^\lambda\leq C_{18}(1+\frac{1}{R^\lambda}) , \nonumber
    \end{eqnarray}
    where $C_{16},C_{17},C_{18}$ are constants depending on $m,k,k_1,\delta_0,\lambda$. 
    \par 
    From the above discussion, we conclude that 
    \begin{eqnarray}
       (\phi \mathcal{F})(x,t)\leq  C_{19}(1+\frac{1}{R^\lambda}) 
    \end{eqnarray}
   on $B_p(2R)\times [0,1]$, where $C_{19}$ is a constant depending on $m,k,k_1,\delta,\delta_0,\lambda$. It follows from (3.45) that 
    \begin{eqnarray}
        \mathcal{F}(x,t)\leq  C_{19}(1+\frac{1}{R^\lambda})  
    \end{eqnarray}
    on $B_p(R)\times [0,1]$. 
    {\qed }
 ~\\
 
    The remaining part of this section is devoted to the case of $t\geq 1$, in which we will consider the auxiliary function
     \begin{eqnarray}
      \mathcal{G}=t(|\nabla_b f|^2+(1+f_0^2)^\lambda-\delta f_t). \nonumber
    \end{eqnarray}
    The argument for this case is almost the same as that for $0<t\leq 1$. Note that at $t=1$, we have
    \begin{align*}
        \mathcal{F}(\cdot,1)=\mathcal{G}(\cdot,1).
    \end{align*}
     
      ~\\
    $\mathbf{Proposition\ 3.4}$ Let $(M^{2m+1},HM,J,\theta)$ be a complete noncompact pseudo-Hermitian manifold with
    \begin{eqnarray}
        Ric_b +2(m-2)Tor_b\geq -k\ \  and\ \   |A|,|\nabla_b A|\leq k_1.\nonumber
     \end{eqnarray}
     and $u$ be a positive solution of the CR heat equation
    \begin{eqnarray}
        \frac{\partial u}{\partial t}=\Delta_b u  \nonumber
     \end{eqnarray}
    on $B_p(2R)\times [1,T]$ with $R\geq 1$ and $T>1$. Then for any constant $\frac{1}{2}<\lambda<\frac{2}{3}$ and any constant $\delta>1+\frac{4}{m\lambda(2\lambda-1)}$, there exists a constant $C_4^{'}$ depending on $m,k,k_1,\lambda,\delta$, such that 
    \begin{eqnarray}
      \frac{|\nabla_b u|^2}{u^2}+(1+\frac{u^2_0}{u^2})^\lambda-\delta \frac{u_t}{u}\leq C_4^{'}(1+\frac{1}{t}+\frac{1}{R^\lambda})
    \end{eqnarray}
    on $B_p(R)\times [1,T]$.
    
    ~\\
     $\mathbf{Proof}$ 
     Let $(x_1,t_1)$ be the maximum point of $\phi \mathcal{G}$ on $M\times [1,T]$. We may assume that $\phi \mathcal{G}$ is positive at $(x_1,t_1)$ and $t_1>1$, otherwise the result follows trivially. Evaluating the inequality in Lemma 3.2 for $\phi \mathcal{G}$ at $(x_1,t_1)$  gives the following
     \begin{eqnarray}
       0&\geq& -\frac{3C_1}{R}\phi^{6\lambda-3}\mathcal{G}
         -\frac{\phi \mathcal{G}}{t_1}+2\langle \nabla_b f,\nabla_b \phi\rangle \mathcal{G}
         \nonumber\\
         &&+\phi t_1\{\frac{1}{m}(\Delta_b f)^2+4mf^2_0 
        \nonumber\\
        &&-\frac{8}{\lambda(2\lambda-1)}(1+f_0^2)^{1-\lambda}|\nabla_b f |^2 -2(k+k_1)|\nabla_b f|^2
        \nonumber\\
        &&-4k_1\lambda (1+f^2_0)^{\lambda-1}|f_0|\cdot|\nabla_b  f|^2 
        \nonumber\\
        &&-2(k_1+k_1^2)\lambda^2(1+f_0^2)^{2\lambda-2}|f_0^2|\}.
     \end{eqnarray}
     Multiplying (3.48) by $\phi t_1$ and using Cauchy-Schwarz inequality, we have 
     \begin{eqnarray}
       0
        &\geq& -\frac{3C_1}{R}\phi^{6\lambda-2}t_1 \mathcal{G}
         -\phi^2 \mathcal{G}-\frac{2C_1}{R}|\nabla_b f|\cdot \phi^{3\lambda} t_1 \mathcal{G} 
         \nonumber\\
         &&+\phi^2 t^2_1\{ \frac{1}{m}(|\nabla_b f|^2-f_t)^2+4mf^2_0 
        -2\epsilon |\nabla_b f|^4-2(k+k_1)|\nabla_b  f|^2
        \nonumber\\
        &&-\frac{16}{\epsilon\lambda^2(2\lambda-1)^2}(1+f_0^2)^{2-2\lambda}
        \nonumber\\
        &&-(2k_1+2k_1^2+4\epsilon^{-1}  k_1^2)\lambda^2(1+f_0^2)^{2\lambda-2}f_0^2\}.
     \end{eqnarray}
     \par 
     The proof of Proposition 3.4 is almost same as that for Proposition 3.3. In following, we only show some necessary modifications when we try to control $\mathcal{G}$. The discussion is similarly divided into the following cases.
     
     ~\\
     $\mathbf{Case\ \Tilde{A}}$ $f_t< 0$ at the maximum point $(x_1,t_1)$ of $\phi \mathcal{G}$. Corresponding to (3.23), we get the following inequality:
     \begin{eqnarray}
       0&\geq&-\frac{3C_1}{R}\phi^{6\lambda-2}t_1 \mathcal{G}
       -\phi^2 \mathcal{G}+2m\phi^2t_1^2f_0^2+\frac{\phi^2 t_1^2}{2m}(|\nabla_b f|^2-f_t)^2
         \nonumber\\
         &&+\phi^2 t^2_1\{ (\frac{1}{2m}-2\epsilon)|\nabla_b f|^4
        -2k|\nabla_b  f|^2-2k_1|\nabla_b  f|^2\}
        \nonumber\\
        &&+\phi^2 t^2_1\{2mf_0^2-\frac{16}{\epsilon\lambda^2(2\lambda-1)^2}(1+f_0^2)^{2-2\lambda}
        \nonumber\\
        &&-(2k_1+2k_1^2+4\epsilon^{-1}  k_1^2)\lambda^2(1+f_0^2)^{2\lambda-2}f_0^2\}
        \nonumber\\
        &&+\phi^2 t^2_1 \frac{1}{2m}f_t^2-\frac{2C_1}{R}|\nabla_b f|\cdot \phi^{3\lambda}t_1 \mathcal{G}.
     \end{eqnarray}
     
     ~\\
     $\mathbf{(\Tilde{A}1)}$ Suppose $\delta_0|\nabla_b f|^2-\delta f_t\geq (1+f^2_0)^\lambda$. Let $\epsilon<\frac{1}{4m}$. Similar to getting
    (3.26) from (3.25), we obtain from (3.50) that
     \begin{eqnarray}
       0&\geq&-\frac{3C_1}{R}\phi^{6\lambda-2}t_1 \mathcal{G}
         -\phi^2 \mathcal{G}+\frac{\phi^2 t_1^2}{2m}(|\nabla_b f|^2-f_t)^2
         \nonumber\\
         &&-C_{20}t^2_1\{1+\frac{1}{R^2}+\frac{1}{R^4}\},
    \end{eqnarray}
    where $C_{20}$ is a constant depending on $m,k,k_1,\delta,\delta_0,\lambda$. Let $x=\phi (\delta_0 |\nabla_b  f|^2-\delta f_t)(x_1,t_1)$, then 
    \begin{eqnarray}
      \phi \mathcal{G}\leq 2t_1 x,\  \phi(|\nabla_b f|^2-f_t)\geq \frac{1}{\delta} x.\nonumber
    \end{eqnarray}
    Consequently we have
    \begin{eqnarray}
      0\geq \frac{t_1^2}{2m\delta^2}x^2-2t_1x-\frac{6C_1}{R}t_1^2 x-C_{20}t_1^2(1+\frac{1}{R^2}+\frac{1}{R^4}).
    \end{eqnarray}
    This implies that 
    \begin{eqnarray}
      t_1 x \leq C_{21}(1+t_1+\frac{t_1}{R}) ,\nonumber
    \end{eqnarray}
    and 
    \begin{eqnarray}
      \phi \mathcal{G}\leq 2t_1 x \leq 2C_{21}(1+t_1+\frac{t_1}{R}) ,
    \end{eqnarray}
     where $C_{21}$ is a constant depending on $m,k,k_1,\lambda,\delta_0,\delta$.
     
     ~\\
     $\mathbf{(\Tilde{A}2)}$ Suppose $\delta_0|\nabla_b f|^2-\delta f_t\leq (1+f^2_0)^\lambda$. Then
      \begin{eqnarray}
      \phi \mathcal{G}
      \leq 2\phi t_1(1+f_0^2)^\lambda. \nonumber
    \end{eqnarray}
    Similar to getting (3.29) from (3.23), we have
     \begin{eqnarray}
       0&\geq&-\frac{3C_1}{R}\phi^{6\lambda-2}t_1 \mathcal{G}
         -\phi^2 \mathcal{G}+2m\phi^2t_1^2f_0^2
         \nonumber\\
         &&+ t^2_1\{ (\frac{1}{2m}-2\epsilon)(\phi|\nabla_b  f|^2)^2
        -2(k+k_1)\phi|\nabla_b  f|^2 \}
        \nonumber\\
        &&+ t^2_1\{2m(\phi|f_0|)^2-\frac{16}{\epsilon\lambda^2(2\lambda-1)^2}-\frac{16}{\epsilon\lambda^2(2\lambda-1)^2}(\phi |f_0|)^{4-4\lambda}
        \nonumber\\
        &&-(2k_1+2k_1^2+4\epsilon^{-1}  k_1^2)\lambda^2(\phi|f_0|)^{4\lambda-2}\}
        \nonumber\\
        &&-\frac{4C_1}{R\sqrt{\delta_0}} [\phi^2(1+f_0^2)]^{\frac{3}{2}\lambda}\}.\
    \end{eqnarray}
    Choosing $\epsilon<\frac{1}{4m}$ and noting that  $\frac{1}{2}<\lambda <\frac{2}{3}$, (3.54) implies that 
    \begin{eqnarray}
       0&\geq&-\frac{3C_1}{R}\phi^{6\lambda-2}t_1 \mathcal{G}
         -\phi^2 \mathcal{G}+2m\phi^2t_1^2f_0^2-C_{22}t_1^2(1+\frac{1}{R}),
         \nonumber
    \end{eqnarray}
    where $C_{22}$ is a constant depending on $m,k,k_1,\delta_0,\lambda$. Set $y=\phi |f_0|$. Then
    \begin{eqnarray}
       0\geq2mt_1^2 y^2-(2t_1+\frac{6C_1t_1^2}{R}) y^{2\lambda}- C_{23}t_1^2(1+\frac{1}{R}),
         \nonumber
    \end{eqnarray}
    where $C_{23}$ is a constant depending on $m,k,k_1,\delta_0,\lambda$. Since $t_1> 1$, we have the following inequality
    \begin{eqnarray}
       0\geq2m y^2-(2+\frac{6C_1}{R}) y^{2\lambda}- C_{23}(1+\frac{1}{R}+\frac{1}{R^4}),
         \nonumber
    \end{eqnarray}
    which yields that 
    \begin{eqnarray}
       y&\leq& C_{24}(1+\frac{1}{\sqrt{R}}+\frac{1}{R^2}+R^{\frac{1}{2(\lambda-1)}} )
         \nonumber\\
        &\leq& C_{25}(1+\frac{1}{\sqrt{R}}),\nonumber
    \end{eqnarray}
   and thus 
    \begin{eqnarray}
       \phi \mathcal{G}\leq 2\phi t_1(1+f^2_0)^\lambda\leq C_{25}(1+\frac{1}{R^\lambda})t_1 ,
    \end{eqnarray}
   where $C_{24},C_{25}$ are constants depending on $m,k,k_1,\delta_0,\lambda$. 
    
   ~\\
   $\mathbf{Case\ \Tilde{B}}$ $f_t\geq 0$ at maximum point $(x_1,t_1)$ of $\phi \mathcal{G}$.
   
   ~\\
   $\mathbf{(\Tilde{B}1)}$ Suppose $\delta_0|\nabla_b f|^2-\delta f_t\geq 0$. We have
   \begin{eqnarray}
       0 &\geq& -\frac{3C_1}{R}\phi^{6\lambda-2}t_1 \mathcal{G}
         -\phi^2 \mathcal{G}+\frac{\phi^2t_1^2}{m\delta^2}(\delta_0|\nabla_b  f|^2-\delta f_t)^2
         \nonumber\\
         &&+\phi^2 t^2_1\{ 
        (\frac{(\delta-\delta_0)^2}{m\delta^2}-2\epsilon )|\nabla_b  f|^4-2(k+k_1)|\nabla_b f|^2\}
        \nonumber\\
        &&+\phi^2t_1^2\{4mf^2_0 -\frac{16}{\epsilon\lambda^2(2\lambda-1)^2}(1+f_0^2)^{2-2\lambda}t_1^{2-4\lambda} 
        \nonumber\\
        &&-4\epsilon^{-1}k^2_1\lambda^2 (1+f^2_0)^{2\lambda-2}|f_0|^2
        \nonumber\\
        &&-2(k_1+k_1^2)\lambda^2(1+f_0^2)^{2\lambda-2}|f_0^2|\}
        \nonumber\\
        &&-\frac{2C_1}{R}|\nabla_b f|\cdot \phi^{3\lambda} t_1 \mathcal{G} .
     \end{eqnarray}
     
     ~\\
     $\mathbf{(\Tilde{B}1}$-$\mathbf{1)} $ Suppose $\delta_0|\nabla_b f|^2-\delta f_t\geq (1+f^2_0)^\lambda$. Similar to (3.36), we have
     \begin{eqnarray}
       0&\geq&-\frac{3C_1}{R}\phi^{6\lambda-2}t_1 \mathcal{G}
         -\phi^2 \mathcal{G}+\frac{\phi^2 t_1^2}{m\delta^2}(\delta_0|\nabla_b f|^2-\delta f_t)^2
         \nonumber\\
         &&-t^2_1\{ C_{26}+\frac{C_{26}}{R}\},
    \end{eqnarray}
    where $C_{26}$ is a constant depending on $m,k,k_1,\delta,\delta_0,\lambda$. Let $x=\phi (\delta_0 |\nabla_b  f|^2-\delta f_t)(x_1,t_1)$. Then we may get
     \begin{eqnarray}
      t_1 x \leq C_{27}(1+t_1+\frac{t_1}{R}) ,\nonumber
    \end{eqnarray}
    that is,
    \begin{eqnarray}
      \phi \mathcal{G}\leq 2t_1 x \leq 2C_{27}(1+t_1+\frac{t_1}{R}) ,
    \end{eqnarray}
     where $C_{27}$ is a constant depending on $m,k,k_1,\lambda,\delta_0,\delta$.
     
     ~\\
   $\mathbf{(\Tilde{B}1}$-$\mathbf{2)}$  Suppose $\delta_0|\nabla_b f|^2-\delta f_t\leq (1+f^2_0)^\lambda$. Corresponding to (3.40), we get the following
    \begin{eqnarray}
       0&\geq&-\frac{3C_1}{R}\phi^{6\lambda-2}t_1 \mathcal{G}
         -\phi^2 \mathcal{G} +2m\phi^2t_1^2f_0^2-C_{28}(1+\frac{1}{R}),\nonumber
    \end{eqnarray}
     where $C_{28}$ is a constant depending on $m,k,k_1,\delta,\delta_0,\lambda$. Set $y=\phi |f_0|$. Then
     \begin{eqnarray}
       0\geq2mt_1^2 y^2-(2t_1+\frac{6C_1t_1^2}{R}) y^{2\lambda}- C_{29}t_1^2(1+\frac{1}{R}),
    \end{eqnarray}
    where $C_{29}$ is a constant depending on $m,k,k_1,\delta,\delta_0,\lambda$. Since $t_1> 1$, we find that  
    \begin{eqnarray}
       0\geq2m y^2-(2+\frac{6C_1}{R}) y^{2\lambda}- C_{30}(1+\frac{1}{R}),
         \nonumber
    \end{eqnarray}
    which yields that 
    \begin{eqnarray}
       y\leq& C_{31}(1+\frac{1}{\sqrt{R}}),\nonumber
    \end{eqnarray}
    hence 
    \begin{eqnarray}
       \phi \mathcal{G}\leq 2\phi t_1(1+f^2_0)^\lambda\leq C_{32}(1+\frac{1}{R^\lambda})t_1 ,
    \end{eqnarray}
    where $C_{30},C_{31},C_{32}$ are constants depending on $m,k,k_1,\delta,\delta_0,\lambda$. 
     
    ~\\
   $\mathbf{(\Tilde{B}2)}$ Suppose $\delta_0|\nabla_b f|^2-\delta f_t\leq 0$. Assume that $\delta_0>1+\frac{4}{m\lambda(2\lambda-1)}$. Similar to getting (3.44) from (3.42) and (3.43), we obtain from (3.48) that
    \begin{eqnarray}
       0&\geq&-\frac{3C_1}{R}\phi^{6\lambda-2}t_1 \mathcal{G}
         -\phi^2 \mathcal{G}+2m\phi^2t_1^2f_0^2-C_{33}(1+\frac{1}{R}),
         \nonumber
    \end{eqnarray}
    where $C_{33}$ is a constant depending on $m,k,k_1,\delta,\delta_0,\lambda$.
     Set $y=\phi |f_0|$, we get
    \begin{eqnarray}
       0\geq2mt_1^2 y^2-(2t_1+\frac{3C_1t_1^2}{R}) y^{2\lambda}- C_{34}(1+\frac{1}{R})
         \nonumber
    \end{eqnarray}
    which yields that 
    \begin{eqnarray}
     y \leq C_{35}(1+\frac{1}{\sqrt{R}}) ,\nonumber
    \end{eqnarray}
    and thus
    \begin{eqnarray}
       \phi \mathcal{G}\leq \phi t_1(1+f^2_0)^\lambda\leq C_{36}(1+\frac{1}{R^\lambda})t_1 ,
    \end{eqnarray}
    where $C_{34},C_{35},C_{36}$ are constants depending on $m,k,k_1,\delta,\delta_0,\lambda$. 
    \par 
    From the above discussion, we conclude that 
    \begin{eqnarray}
       \phi(x) \mathcal{G}(x,t)\leq  C_{37}(1+t_1+\frac{t_1}{R^\lambda}) ,
    \end{eqnarray}
    on $B_p(2R)\times [1,T]$, where $C_{37}$ is a constant depending on $m,k,k_1,\delta,\delta_0,\lambda$. Consequently
    \begin{eqnarray}
        \mathcal{G}(x,t)\leq  C_{37}(1+t_1+\frac{t_1}{R^\lambda})
    \end{eqnarray}
    on $B_p(R)\times [1,T]$. In particular, we have 
    \begin{eqnarray}
        \mathcal{G}(x,T)\leq  C_{37}(1+t_1+\frac{t_1}{R^\lambda})\leq C_{37}(1+T+\frac{T}{R^\lambda})
    \end{eqnarray}
    on $B_p(R)$. Since $T(>1)$ is arbitrary, this gives (3.47).
    {\qed }
    
    ~\\
     Combining Propositions 3.3 and 3.4,  we may obtain Theorem 1.1. Clearly Theorem 1.2 follows from Theorem 1.1 by letting $R\rightarrow \infty$. We would like to end this section by the following remark.
    
    ~\\
    $\mathbf{Remark\ 3.1.}$ If $M$ is a closed pseudo-Hermitian manifold with the same properties as in Theorem 1.1, we may carry out the above argument, without using the cut-off function, to deduce the following result: Let $u$ be a positive solution of the CR heat equation on $M^{2m+1}$. Then for any constant $\frac{1}{2}<\lambda<\frac{2}{3}$ and any constant $\delta>1+\frac{4}{m\lambda(2\lambda-1)}$, there exists a constant $C$ depending on $m,k,k_1,\lambda,\delta$, such that 
    \begin{align}
        \frac{|\nabla_b u|^2}{u^2}-\delta \frac{u_t}{u}\leq C+\frac{C}{t}
    \end{align}
    on $M\times (0,\infty)$. We should point out that the sub-Laplacian $\Delta_b$ can only be expressed as (1.1) locally, that is,
     \begin{align}
        \Delta_b=\sum\limits_{A=1}^{2m}e_A^2-\sum\limits_{A=1}^{2m}\nabla_{e_A}e_A, 
    \end{align}
    where $\{e_A\}_{A=1}^{2m}$ is the local frame field given in $\S 2$, and $\nabla$ is the Tanaka-Webster connection. In general, one cannot express $ \Delta_b$ as (3.66) by global vector fields. Hence, although the method for the closed case follows essentially from Cao and Yau \cite{CaY}, their result cannot be applied directly to get the estimate (3.65).

    ~\\
    $\mathbf{Remark\ 3.2.}$ Note that $1+\frac{4}{m\lambda(2\lambda-1)}>1+\frac{18}{m}$ for $ \frac{1}{2}<\lambda<\frac{2}{3}$. Let $M^{2m+1}$ be either a closed pseudo-Hermitian manifold or a complete noncompact pseudo-Hermitian manifold with the same properties as in Theorem 1.1. Then we have the following Li-Yau type estimate: Let $u$ be a positive solution of the CR heat equation on $M$. Then for any $\delta>1+\frac{18}{m}$, there exists a constant $C$ depending on $m,k,k_1,\delta$, such that 
    \begin{align}
        \frac{|\nabla_b u|^2}{u^2}-\delta \frac{u_t}{u}\leq C+\frac{C}{t}
    \end{align}
    on $M\times (0,\infty)$.
    \section{Harnack inequality and heat kernel estimates}
    In this section, we derive the CR version of Harnack's inequality for the positive solutions of the CR heat equation and deduce an upper bound for the heat kernel.
    
    ~\\
    $\mathbf{Proof\ of\ Theorem\ 1.3}$ Let $\gamma:[t_1,t_2]\rightarrow M$ be a horizontal curve joining $x$ and $y$, i.e. $\gamma(t_1)=x,\gamma (t_2)=y$. Define a map $\eta:[t_1,t_2]\rightarrow M\times [t_1,t_2]$ by 
     \begin{eqnarray}
      \eta(t)=(\gamma(t),t). \nonumber
    \end{eqnarray} 
    Let $f=\mathrm{ln}\ u$ with $u$ being a positive solution of the CR heat equation. Integrating $\frac{d}{dt}f $ along $\eta$, we get
    \begin{eqnarray}
      f(y,t_2)-f(x,t_1)=\int_{t_1}^{t_2}\frac{d}{dt}f dt=\int_{t_1}^{t_2}(\langle \Dot{\gamma},\nabla_b f\rangle +f_t)dt. \nonumber
    \end{eqnarray} 
    Applying Theorem 1.2 gives 
    \begin{eqnarray}
      f(y,t_2)-f(x,t_1)&\geq& \int_{t_1}^{t_2}(\langle \Dot{\gamma},\nabla_b f\rangle+\frac{1}{\delta}|\nabla_bf|^2-\frac{C}{\delta}-\frac{C}{\delta t} )dt, \nonumber \\
      &\geq & -\int_{t_1}^{t_2}(\frac{\delta |\Dot{\gamma }|^2}{4}+\frac{C}{\delta}+\frac{C}{\delta t} )dt.\nonumber
    \end{eqnarray}
   Choosing a curve $\gamma$ with $|\Dot{\gamma}|=\frac{d_{cc}(x,y)}{t_2-t_1}$, we obtain 
    \begin{eqnarray}
      \mathrm{ln}\frac{u(y,t_2)}{u(x.t_1)}\geq -\frac{C}{\delta}(t_2-t_1)-\frac{C}{\delta}\mathrm{ln}\frac{t_2}{t_1}-\frac{\delta d_{cc}^2(x,y)}{4(t_2-t_1)}. \nonumber
    \end{eqnarray}
    Taking exponentials of the above inequality, we can complete the proof.
    {\qed}
    \par
    A mean value type inequality follows immediately from Theorem 1.3.
    
    ~\\
    $\mathbf{Corollary\ 4.1}$ Let $(M^{2m+1},HM,J,\theta)$ be a complete noncompact pseudo-Hermitian manifold with
    \begin{eqnarray}
        Ric_b+2(m-2)Tor_b\geq -k,\ and\   |A|,|\nabla_b A|\leq k_1,\nonumber
     \end{eqnarray}
     and $u$ be a positive solution of the heat equation
    \begin{eqnarray}
        \frac{\partial u}{\partial t}=\Delta_b u  \nonumber
     \end{eqnarray}
    on $M\times (0,\infty)$. Then for any constant $\frac{1}{2}<\lambda<\frac{2}{3}$ and any constant $\delta>1+\frac{4}{m\lambda(2\lambda-1)}$, there exists a constant $C$ which is given by Theorem 1.2 such that for any $0<t_1<t_2$ and $x\in M$, we have
    \begin{eqnarray}
      u(x,t_1)&\leq& [Vol(B_{cc}(x,r))]^{-\frac{1}{2}}\left ( \int_{B_{cc}(x,r)}u^2(y,t_2)dy\right) ^{\frac{1}{2}}(\frac{t_2}{t_1})^{\frac{C}{\delta}} \nonumber\\ 
      && \cdot exp(\frac{C}{\delta}(t_2-t_1)+\frac{\delta r^2}{4(t_2-t_1)} ).
    \end{eqnarray} 
    \par
    Our next goal of this section is to derive an upper estimate for the heat kernel. For any $x,y \in M$ and $t>0$, let us set $\rho(x,y,t)=\frac{1}{2t}d^2_{cc}(x,y)$. It is known that 
    \begin{eqnarray}
       |\nabla_b d_{cc}|^2\leq 1  \nonumber
     \end{eqnarray}
    in the weak sense, where $ \nabla_b d_{cc}$ denotes the horizontal gradient of $d_{cc}$ with respect to either $x$ or $y$ (\cite{JX}). Define $g(x,y,t)=-\rho(x,y,(1+2\alpha)T-t)$, where $\alpha, T$ are constants to be determined later. A direct computation shows that
    \begin{eqnarray}
       \frac{1}{2}|\nabla_b g|^2+g_t \leq 0. 
     \end{eqnarray}
    
    ~\\ 
    $\mathbf{Lemma\ 4.2}$ Let $(M^{2m+1},HM,J,\theta)$ be a complete pseudo-Hermitian manifold. Suppose $H(x,y,t)$ is the heat kernel of (1.4). Let
     \begin{eqnarray}
       \mathcal{F}_x(y,t)=\int_{S_1}H(y,z,t)H(x,z,T)dz \nonumber
     \end{eqnarray}
    for any fixed $x\in M$ and any nonempty subset $S_1\subset M$. Then for any $0\leq t\leq s<(1+2\alpha)T$ and any nonempty subset $S_2\subset M $, we have 
    \begin{eqnarray}
       \int_{S_2}\mathcal{F}_x^2(z,s)dz &\leq & \int_{S_1}H^2(x,z,T)dz\sup\limits_{z\in S_1}exp(-2\rho(x,z,(1+2\alpha)T)) \nonumber \\
       && \cdot \sup\limits_{z\in S_2} exp (2\rho(x,z,(1+2\alpha)T-s)) .
     \end{eqnarray}
    
    ~\\
    $\mathbf{Proof}$  The proof is almost the same as in \cite{LY} by using (4.2).
    {\qed}
     
    ~\\
    $\mathbf{Proof\ of\ Theorem\ 1.4}$ Let $\mathcal{F}_x$ be defined as in Lemma 4.2, and let $S_1=B_{cc}(y,\sqrt{t}),S_2=B_{cc}(x,\sqrt{t})$. Applying Corollary 4.1 to $\mathcal{F}_x$ with $t_1=T$ and $t_2=(1+\alpha T)$, we have 
    \begin{eqnarray}
        &&(\int_{B_{cc}(y,\sqrt{t})}H^2(x,z,T)dz)^2=\mathcal{F}_x^2(x,T)
        \nonumber\\
        &\leq& [Vol(B_{cc}(x,\sqrt{t}))]^{-1}\int_{B_{cc}(x,\sqrt{t})}\mathcal{F}_x^2(z,(1+\alpha)T)dz\cdot (1+\alpha )^{\frac{2C}{\delta}}
        \nonumber\\
        && \cdot exp( \frac{2C}{\delta}\alpha T+\frac{\delta t}{2\alpha T}). 
     \end{eqnarray}
    By Lemma 4.2, we deduce from (4.4) that 
    \begin{eqnarray}
        \mathcal{F}_x^2(x,T)&\leq& [Vol(B_{cc}(x,\sqrt{t}))]^{-1}\int_{B_{cc}(y,\sqrt{t})}H^2(x,z,T)dz\cdot (1+\alpha )^{\frac{2C}{\delta}} \nonumber\\
        &&\cdot exp( \frac{2C}{\delta}\alpha T+\frac{(\delta+1) t}{2\alpha T}-\inf\limits_{z\in B_{cc}(y,\sqrt{t})}\rho(x,z,(1+2\alpha)T)). \nonumber
     \end{eqnarray}
    Hence 
    \begin{align*}
        \int_{B_{cc}(y,\sqrt{t})}H^2(x,z,T)dz\leq& [Vol(B_{cc}(x,\sqrt{t}))]^{-1}(1+\alpha )^{\frac{2C}{\delta}}exp( \frac{2C}{\delta}\alpha T+\frac{(\delta+1) t}{2\alpha T}) \nonumber\\
        &\cdot exp( -\inf\limits_{z\in B_{cc}(y,\sqrt{t})}\rho(x,z,(1+2\alpha)T)). \nonumber
    \end{align*}
    Applying Corollary 4.1 once again and letting $T=(1+\alpha)t$, we have
    \begin{eqnarray}
        &&H^2(x,y,t)
        \nonumber\\
        &\leq&[Vol(B_{cc}(y,\sqrt{t}))]^{-1}\int_{B_{cc}(y,\sqrt{t})}H^2(x,y,T)dz(1+\alpha )^{\frac{2C}{\delta}}
        \ exp( \frac{2C}{\delta}\alpha t+\frac{\delta }{2\alpha }) \nonumber\\
        &\leq&[Vol(B_{cc}(y,\sqrt{t}))]^{-1}[Vol(B_{cc}(x,\sqrt{t}))]^{-1}(1+\alpha )^{\frac{4C}{\delta}} \nonumber\\
        &&\cdot  exp\left( \frac{4C}{\delta}(\alpha+2)\alpha t +\frac{(\alpha+2)\delta+1 }{2\alpha(1+\alpha) }-\inf\limits_{z\in B_{cc}(y,\sqrt{t})}\rho(x,z,(1+2\alpha)T)\right). \nonumber
     \end{eqnarray}
    Similar to \cite{CaY}, we set 
     $4(1+\alpha)^2(1+2\alpha)=4+\epsilon$, then 
    \begin{eqnarray}
       \inf\limits_{z\in B_{cc}(y,\sqrt{t})}\rho(x,z,(1+2\alpha)T)\geq \frac{2d_{cc}^2(x,y)}{(4+\epsilon)t}-\frac{2(1+\alpha)}{(4+\epsilon)\alpha}.
     \end{eqnarray}
     This completes the proof.
     {\qed}
     
     ~\\
     $\mathbf{Remark\ 4.1.}$  A similar result for subelliptic operators on closed manifolds was given in \cite{CaY}.

\bibliographystyle{siam}
\bibliography{ref}

~\\
  Yuxin Dong
  \\
  $School\ of\ Mathematical\ Science $
  \\
  and 
  \\
  $Laboratory\ of\ Mathematics\ for\ Nonlinear\ Science $
\\
  $ Fudan\ University$
\\
   $Shanghai\ 200433 ,$    $P.R.\ China $
   \\
   $yxdong@fudan.edu.cn$

   ~\\
   Yibin Ren
   \\
   $College\ of Mathematics$
   \\
   $Physics\ and\ Information\ Engineering$
   \\
   $Zhejiang\ Normal\ University$
   \\
   $Jinhua\ 321004,$ $P.R.\ China $
   \\
  $ allenry@outlook.com$

   ~\\
  Biqiang Zhao
  \\
  $Shanghai\ Center\ For$
  \\
  $Mathematical\ Science $
\\
  $ Fudan\ University$
\\
   $Shanghai\ 200433 ,$ $P.R.\ China $
   \\
   $17110840003@fudan.edu.cn$

\end{document}